\documentclass[final]{siamltex1213}
\usepackage{layout}
\usepackage{amsmath,amsfonts,amssymb}
\usepackage{mathtools}
\usepackage{xcolor}
\usepackage{hyperref}
\usepackage{algorithm}
\usepackage{algpseudocode}
\usepackage{subcaption}
\usepackage{multirow}
\usepackage{pgfplots}

\title{Numerically stable variants of the communication-hiding pipelined Conjugate Gradients algorithm for the parallel solution of large scale symmetric linear systems}
\author{S. Cools$^*$, W. Vanroose\thanks{Applied Mathematics Group, Department of Mathematics and Computer Science, University of Antwerp, Middelheimlaan 1, 2020 Antwerp, Belgium. Contact: siegfried.cools@uantwerp.be}} 

\begin{document}

\maketitle

\begin{abstract}
By reducing the number of global synchronization bottlenecks per iteration and hiding communication behind useful computational work, pipelined Krylov subspace methods achieve significantly improved parallel scalability on present-day HPC hardware. However, this typically comes at the cost of a reduced maximal attainable accuracy.
This paper presents and compares several stabilized versions of the communication-hiding pipelined Conjugate Gradients method. The main novel contribution of this work is the reformulation of the multi-term recurrence pipelined CG algorithm by introducing shifts in the recursions for specific auxiliary variables. 
These shifts reduce the amplification of local rounding errors on the residual. The stability analysis presented in this work provides a rigorous method for selection of the optimal shift value in practice. It is shown that, 
given a proper choice for the shift parameter, the resulting shifted pipelined CG algorithm restores the attainable accuracy and displays nearly identical robustness to local rounding error propagation compared to classical CG. 
Numerical results on a variety of SPD benchmark problems compare different stabilization techniques for the pipelined CG algorithm, showing that the shifted pipelined CG algorithm is able to attain a high accuracy while displaying excellent parallel performance.
\end{abstract}


\section{Introduction} 

Both on the academic and industrial level, Krylov subspace methods 
\cite{greenbaum1997iterative,liesen2012krylov,meurant1999computer,saad2003iterative,van2003iterative} 
are well-known as efficient solution methods for large scale linear systems in high-performance computing. 
These iterative algorithms constructs a sequence of approximate solutions $\{x_i\}_i$ with $x_i \in x_0 + \mathcal{K}_i(A,r_0)$ to the system $Ax =b$ , 
where $r_0 = b-Ax_0$ is the initial residual and the $i$-th Krylov subspace is 
\begin{equation}
	\mathcal{K}_i(A,r_0) := \text{span}\{r_0,Ar_0,A^2r_0,\ldots,A^{i-1}r_0\}.
\end{equation}
The Conjugate Gradient (CG) method \cite{hestenes1952methods}, Alg.\,\ref{algo::pcg}, which allows for 
the solution of linear systems with symmetric positive definite (SPD) matrices $A$,
is generally considered as the first Krylov subspace method. 
Driven by the transition of hardware towards the exascale regime, research on the scalability of Krylov subspace methods
on massively parallel architectures has recently gained increasing attention 
\cite{dongarra2011international,dongarra1998numerical,dongarra2013toward, dongarra2015hpcg}.
Since for many applications the system matrix is sparse and thus inexpensive to apply, the main bottleneck for efficient parallel execution is typically not this sparse 
matrix-vector product (\textsc{spmv}), but the communication overhead caused by global reductions in dot-product computations. 

Over the last decades significant efforts have been made to reduce and/or eliminate the synchronization bottleneck in Krylov subspace methods.
The earliest papers on synchronization reduction date back to the late 1980's and 1990's, 
see \cite{barrett1994templates,d1992reducing,de1991parallel,demmel1993parallel,erhel1995parallel,strakovs1987effectivity}.
A reduction of the number of global communication points was also introduced by the so-called $s$-step methods by Chronopoulos 
et al.~\cite{chronopoulos1989s,chronopoulos2010block,chronopoulos1996parallel} and Carson et al.~\cite{carson2014residual,carson2013avoiding}. 
In addition to communication avoiding methods, research on hiding global communication by overlapping communication with computations can be found in the literature, see \cite{de1995reducing,demmel1993parallel,ghysels2013hiding,ghysels2014hiding}. 

The pipelined CG (p-CG) method proposed in \cite{ghysels2014hiding} aims at hiding global synchronization latency 
by overlapping the global communication phase by the \textsc{spmv}, which requires only local communication. In this way, idle core time 
is minimized by performing useful computations simultaneously to the time-consuming synchronization phase, cf.~\cite{eller2015non}.
The reorganization of the algorithm required to achieve the overlap 
introduces several additional \textsc{axpy} ($y \leftarrow \alpha x +y$) operations to recursively compute auxiliary variables. 
Since vector operations such as an \textsc{axpy} are computed locally, they do not require communication between nodes. 
Thus, the addition of extra recurrences has no impact on the communication flow of the algorithm. 

However, extra recurrences may influence the numerical stability of the algorithm.
Whereas in exact arithmetic the pipelined CG algorithm is equivalent to classic CG, in finite precision 
each of the additional recurrences introduce local rounding errors. 
As analyzed in \cite{cools2016analysis}, the propagation of local rounding errors is more severe for pipelined CG 
compared to classic CG, and can have a detrimental effect on the attainable accuracy of the final iterative solution.
Note that, apart from the reduced attainable accuracy, other rounding error effects may affect the convergence of the multi-term recurrence pipelined CG method.
For example, in some applications a delay of convergence due to loss of the Lanczos basis vector orthogonality in finite precision arithmetic may be observed, see \cite{carson2016numerical,gergelits2014composite,greenbaum1992predicting,strakovs2002error}. 

\begin{algorithm}[t]
  \caption{Classic preconditioned CG}
  \label{algo::pcg}
  \begin{algorithmic}[1]
  	\Procedure{cg}{$A$, $M^{-1}$, $b$, $x_0$}
    \State $r_0 := b - Ax_0$; $u_0:= M^{-1} r_0$; $p_0 := u_0$  \label{al:residual}
    \For{$i = 0, \dots$}
    \State $s_i := Ap_{i}$ \label{al:matvec}
    \State $\alpha_{i} := \left( r_i, u_i \right) / \left( s_i, p_i \right)$ \label{al:alpha}
    \State $x_{i+1} := x_i + \alpha_{i} p_i$
    \State $r_{i+1} := r_i - \alpha_{i} s_i$
    \State $u_{i+1} := M^{-1} r_{i+1}$
    \State $\beta_{i+1} := \left( r_{i+1}, u_{i+1} \right) / \left( r_i, u_i \right)$
    \State $p_{i+1} := u_{i+1} + \beta_{i+1} p_i$ 
    \EndFor
    \EndProcedure
  \end{algorithmic}
\end{algorithm}

\section{The shifted pipelined CG algorithm} 

\subsection{Standard pipelined Conjugate Gradients}

\begin{algorithm}[t]
  \caption{Preconditioned pipelined CG}
  \label{algo::ppipe-cg}
  \begin{algorithmic}[1]
  	\Procedure{p-cg}{$A$, $M^{-1}$, $b$, $x_0$}
    \State $r_0 := b - Ax_0$; $u_0:= M^{-1} r_0$; $w_0 := Au_0$
    \For{$i = 0,\dots$}
    \State $\gamma_i :=(r_i,u_i)$
    \State $\delta := (w_i,u_i)$
    \State $m_i := M^{-1} w_i$
    \State $n_i := A m_i$
    \If{$i>0$}
    \State $\beta_i := \gamma_i/\gamma_{i-1}$; $\alpha_i := (\delta/\gamma_i - \beta_i/\alpha_{i-1})^{-1}$
    \Else
    \State $\beta_i := 0$; $\alpha_i := \gamma_i/\delta$
    \EndIf
    \State $z_i := n_i + \beta_i z_{i-1}$
    \State $q_i := m_i + \beta_i q_{i-1}$
    \State $s_i := w_i + \beta_i s_{i-1}$
		\State $t_i := r_i + \beta_i t_{i-1}$
    \State $p_i := u_i + \beta_i p_{i-1}$
    \State $x_{i+1} := x_i + \alpha_i p_i$
    \State $r_{i+1} := r_i - \alpha_i s_i$
    \State $u_{i+1} := u_i - \alpha_i q_i$
    \State $w_{i+1} := w_i - \alpha_i z_i$
    \EndFor
    \EndProcedure
  \end{algorithmic}
\end{algorithm}

The communication hiding pipelined Conjugate Gradient algorithm for solving the SPD system $Ax = b$ is shown in Alg.\,\ref{algo::ppipe-cg}. Here $A$ is assumed to be a real-valued symmetric positive definite $n$-by-$n$ matrix, i.e., $A \in \mathbb{R}^{n \times n}$, and $b \in \mathbb{R}^{n \times 1}$ is the right-hand side vector. 

Only one \textsc{spmv} is performed in each iteration of the (pipelined) CG algorithm. All other auxiliary vectors are defined recursively to reduce computational overhead. Apart from the current solution vector $x_i$, the (unpreconditioned) search direction $t_i$ and the residual $r_i = b - A x_i$, the following auxiliary variables are introduced in the pipelined CG algorithm:
\begin{align}
	u_i &:= M^{-1} r_i , &
	w_i &:= A u_i , &
	m_i &:= M^{-1} w_i ,&
	n_i &:= A m_i , \notag \\
	p_i &:= M^{-1} t_i , &
	s_i &:= A p_i , &
	q_i &:= M^{-1} s_i , &
	z_i &:= A q_i ,\label{eq:pipe_cg_aux}
\end{align}
With the exception of $m_i$ and $n_i$ these variables are defined recursively in Alg.\,\ref{algo::ppipe-cg}. We refer to the original paper on pipelined CG \cite{ghysels2014hiding} and our own contribution \cite{cools2017communication} for details on the derivation of the corresponding recurrences.

The operator $M^{-1}$ denotes the (left or right) preconditioner, which may either be defined explicitly as $M^{-1} \in \mathbb{R}^{n \times n}$ or given as a general vector operator. 
In the unpreconditioned case the operator $M^{-1}$ is the identity matrix $I \in \mathbb{R}^{n \times n}$, and the definitions for the auxiliary variables $u_i$, $m_i$, $p_i$ and $q_i$ in Eq.\,\eqref{eq:pipe_cg_aux} are redundant.

\subsection{Derivation of shifted pipelined Conjugate Gradients} \label{sec:derivation}

Similar to the construction of pipelined CG in \cite{ghysels2014hiding}, we derive the shifted pipelined CG algorithm from the basic recurrences for the search direction $t_i$, solution $x_i$ and residual $r_i$. The following auxiliary variables are introduced:
\begin{align}
	u_i &:= M^{-1} r_i , &
	w_i &:= \left(AM^{-1}-\sigma I\right) r_i = Au_i -\sigma r_i, \notag \\
	m_i &:= M^{-1} w_i ,&
	n_i &:= A m_i , \notag \\
	p_i &:= M^{-1} t_i , &
	s_i &:= \left(AM^{-1}-\sigma I\right) \, t_i = Ap_i -\sigma t_i, \notag \\
	q_i &:= M^{-1} s_i , &
	z_i &:= A q_i , \label{eq:pipe_cg_shift_aux}
\end{align}
The scalar shift parameter $\sigma \in \mathbb{R}_0^+$ is assumed to be strictly positive. The matrix $I \in \mathbb{R}^{n \times n}$ denotes the unit matrix. 
In the unpreconditioned case, the definitions for the shifted auxiliary variables $w_i$ and $s_i$ 
reduce to $w_i = \left(A-\sigma I \right) r_i$ and $s_i = \left(A-\sigma I \right) t_i$. 

We start from the recurrence for the unpreconditioned search direction
\begin{equation}
	t_i = r_i + \beta_i t_{i-1} , \label{eq:hatp} 
\end{equation}
and the preconditioned version, which is derived by multiplying Eq.\,\eqref{eq:hatp} by $M^{-1}$;
\begin{equation}
	p_i = u_i + \beta_i p_{i-1} . 
\end{equation}
Multiplying Eq.\,\eqref{eq:hatp} by the shifted operator $\left(A M^{-1} - \sigma I \right)$ on both sides, we obtain
\begin{equation}
s_i = \left(A M^{-1} - \sigma I \right) \left( r_i + \beta_i t_{i-1} \right) = w_i + \beta_i s_{i-1}.  \label{eq:s} 
\end{equation}
The recursion for $s_i$ is hence identical to the non-shifted case, see Alg.\,\ref{algo::ppipe-cg}, line 15.

Combining the definition of the residual $r_i = b - A x_i$, the recurrence for the solution $x_{i+1} = x_i + \alpha_i p_i$, and the definition of the variable $s_i = Ap_i - \sigma t_i$, we derive the following recurrence for the residual:
\begin{equation}
	r_{i+1} = b - A x_{i+1} = b - A x_i - \alpha_i A p_i = r_i - \alpha_i s_i - \alpha_i \sigma t_i. \label{eq:r} 
\end{equation}
By multiplying the above on both sides by $M^{-1}$, we obtain the recurrence for the preconditioned residual $u_{i+1} = M^{-1} r_{i+1}$;
\begin{equation}
	u_{i+1} = M^{-1} r_i - \alpha_i M^{-1} s_i - \alpha_i \sigma M^{-1}t_i = u_i - \alpha_i q_i - \alpha_i \sigma p_i. \label{eq:u} 
\end{equation}
In the recurrences for $r_i$ and $u_i$, Eq.\,\eqref{eq:r} and Eq.\,\eqref{eq:u}, a term that depends on the shift parameter $\sigma$ is subtracted to compensate for the introduction of the shift in $s_i$.

\begin{algorithm}[t]
  \caption{Shifted preconditioned pipelined CG}
  \label{algo::ppipe-cg-shift}
  \begin{algorithmic}[1]
  	\Procedure{p-cg-sh}{$A$, $M^{-1}$, $b$, $x_0$, $\sigma$}
    \State $r_0 := b - Ax_0$; $u_0:= M^{-1} r_0$; $w_0 := Au_0 - \sigma r_0$
    \For{$i = 0,\dots$}
    \State $\gamma_i :=(r_i,u_i)$
    \State $\delta := (w_i+\sigma r_i,u_i)$
    \State $m_i := M^{-1} w_i$
    \State $n_i := A m_i$
    \If{$i>0$}
    \State $\beta_i := \gamma_i/\gamma_{i-1}$; $\alpha_i := (\delta/\gamma_i - \beta_i/\alpha_{i-1})^{-1}$
    \Else
    \State $\beta_i := 0$; $\alpha_i := \gamma_i/\delta$
    \EndIf
    \State $z_i := n_i + \beta_i z_{i-1}$
    \State $q_i := m_i + \beta_i q_{i-1}$
    \State $s_i := w_i + \beta_i s_{i-1}$
		\State $t_i := r_i + \beta_i t_{i-1}$
    \State $p_i := u_i + \beta_i p_{i-1}$
    \State $x_{i+1} := x_i + \alpha_i p_i$
    \State $r_{i+1} := r_i - \alpha_i s_i - \alpha_i \sigma t_i$
    \State $u_{i+1} := u_i - \alpha_i q_i - \alpha_i \sigma p_i$
    \State $w_{i+1} := w_i - \alpha_i z_i$
    \EndFor
    \EndProcedure
  \end{algorithmic}
\end{algorithm}

Applying the shifted operator $\left(A M^{-1} - \sigma I\right)$ to both sides in Eq.\,\eqref{eq:r} yields
\begin{equation}
w_{i+1} = w_i - \alpha_i \left(A M^{-1} - \sigma I\right) s_i - \alpha_i \sigma s_i = w_i - \alpha_i z_i. \label{eq:w} 
\end{equation}
To obtain the recurrence for $q_i$, we multiply Eq.\,\eqref{eq:s} by $M^{-1}$, which results in
\begin{equation}
	q_i = M^{-1} w_i + \beta_i M^{-1} s_{i-1} = m_i + \beta_i q_{i-1}.  \label{eq:q} 
\end{equation}
Finally, the recurrence for the variable $z_i$ is found by multiplying Eq.\,\eqref{eq:q} by $A$, i.e.;
\begin{equation}
	z_i = A m_i + \beta_i A q_{i-1} = n_i + \beta_i z_{i-1} . \label{eq:z} 
\end{equation}
The resulting algorithm, denoted as shifted pipelined CG, is summarized in Alg.\,\ref{algo::ppipe-cg-shift}. Alg.\,\ref{algo::ppipe-cg-shift} reverts to standard pipelined CG, Alg.\,\ref{algo::ppipe-cg}, when the shift is set to $\sigma = 0$.

\section{Numerical stability analysis for shifted pipelined CG} \label{sec:rounding}

The pipelined variants, Alg.\,\ref{algo::ppipe-cg} and Alg.\,\ref{algo::ppipe-cg-shift}, are fully equivalent to classic CG, Alg.\,\ref{algo::pcg}, in exact arithmetic. Indeed, in exact arithmetic the recurrences for the residuals (and other auxiliary variables) are identical to the explicit \textsc{spmv}-based formulations, e.g., $r_{i+1} = b - A x_{i+1} = r_i - \alpha_i s_i - \alpha_i \sigma t_i$ holds.
However, when implemented in finite precision arithmetic in practice, local rounding errors contaminate the recurrences, inducing a gap between the explicit and recursive characterizations. In this section we consider a finite precision framework, in which computed variables are denoted by a bar symbol.
The following analysis builds on the work by Greenbaum \cite{greenbaum1997estimating} and Strako{\v{s}} \& Gutknecht \cite{gutknecht2000accuracy}. Strongly related work can be found in \cite{demmel1997applied,greenbaum1989behavior,meurant2006lanczos,paige1971computation,paige1972computational,paige1976error,paige1980accuracy,strakovs2002error,strakovs2005error}.

\subsection{Analysis of local rounding error propagation}

To analyze the propagation of local rounding errors introduced by the recurrences in p-CG-sh, Alg.\,\ref{algo::ppipe-cg-shift}, the following model for floating point arithmetic with machine precision $\epsilon$ is assumed:
\begin{equation}
	\text{fl}(a\pm b) = a(1+\epsilon_1) \pm b(1+\epsilon_2), \quad |\epsilon_1|,|\epsilon_2| \leq \epsilon,
\end{equation}
\vspace{-0.8cm}
\begin{equation}
	\text{fl}(a \text{~op~} b) = (a \text{~op~} b) (1+\epsilon_3), \quad |\epsilon_3| \leq \epsilon, \quad \text{op} = *, /.
\end{equation}
Under this model, and discarding terms involving $\epsilon^2$ or higher powers of $\epsilon$ when terms of order $\epsilon$ are present, the following standard results for operations on an $n$-by-$n$ matrix $A$, $n$-length vectors $v$ and $w$ and a scalar $\alpha$ hold:
\begin{equation}
	\| \alpha v - \text{fl}(\alpha v) \| \leq \| \alpha v \| \, \epsilon =  |\alpha| \, \|v\| \, \epsilon,
\end{equation}
\vspace{-0.8cm}
\begin{equation}
	\| v + w - \text{fl}(v + w) \| \leq (\|v\| + \|w\|) \, \epsilon,
\end{equation}
\vspace{-0.8cm}
\begin{equation}
	| \left( v,w \right) - \text{fl}(\,\left(v,w \right)\,) | \leq n \, \|v\| \, \|w\| \epsilon,
\end{equation}
\vspace{-0.8cm}
\begin{equation}
	\| Av - \text{fl}(Av) \| \leq (\mu\sqrt{n}) \, \|A\| \, \|v\| \, \epsilon ,
\end{equation}
where $\mu$ is the maximum number of nonzeros in any row of $A$. The norm $\|\cdot \|$ denotes the Euclidean 2-norm.

Replacing all recurrences in shifted preconditioned pipelined CG, Alg.\,\ref{algo::ppipe-cg-shift}, by their finite precision equivalents, we obtain:
\begin{align}
	\bar{x}_{i+1} &= \bar{x}_i + \bar{\alpha}_i \bar{p}_i + \delta_i^x, &
	\bar{t}_i &= \bar{r}_i + \bar{\beta}_i \bar{t}_{i-1} + \delta_i^{t} , \notag \\
	\bar{r}_{i+1} &= \bar{r}_i - \bar{\alpha}_i \bar{s}_i -\bar{\alpha}_i \sigma \bar{t}_i+ \delta_i^r,&
	\bar{p}_i 		&= \bar{u}_i + \bar{\beta}_i \bar{p}_{i-1} + \delta_i^p , \notag \\
	\bar{u}_{i+1} &= \bar{u}_i - \bar{\alpha}_i \bar{q}_i - \bar{\alpha}_i \sigma \bar{p}_i + \delta_i^u , &
	\bar{s}_i 		&= \bar{w}_i + \bar{\beta}_i \bar{s}_{i-1} + \delta_i^s, \notag \\
	\bar{w}_{i+1} &= \bar{w}_i - \bar{\alpha}_i \bar{z}_i + \delta_i^w, &	
	\bar{q}_i     &= \bar{m}_i + \bar{\beta}_i \bar{q}_{i-1} + \delta_i^q, \notag \\
	& &
	\bar{z}_i     &= \bar{n}_i + \bar{\beta}_i \bar{z}_{i-1} + \delta_i^z, \label{eq:recur}
\end{align}
where the local rounding errors on each variable are bounded as follows
\begin{eqnarray} 
	\|\delta_i^x\| &\leq& \left( \|\bar{x}_i\| + 2 \, |\bar{\alpha}_i| \, \|\bar{p}_i\| \right) \epsilon , \notag \\
	\|\delta_i^{t}\| &\leq& \left( \|\bar{r}_i\| + 2 \, |\bar{\beta}_i| \, \|\bar{t}_{i-1}\| \right) \epsilon , \notag \\
	\|\delta_i^r\| &\leq& \left( \|\bar{r}_i\| + 3 \, |\bar{\alpha}_i| \, \|\bar{s}_i\| + 4 \, |\bar{\alpha}_i| \, |\sigma| \, \|\bar{t}_i\| \right) \epsilon , \notag \\
	\|\delta_i^p\| &\leq& \left( \|\bar{u}_i\| + 2 \, |\bar{\beta}_i| \, \|\bar{p}_{i-1}\| \right) \epsilon , \notag \\
	\|\delta_i^u\| &\leq& \left( \|\bar{u}_i\| + 3 \, |\bar{\alpha}_i| \, \|\bar{q}_i\| + 4 \, |\bar{\alpha}_i| \, |\sigma| \, \|\bar{p}_i\| \right) \epsilon  , \notag \\
	\|\delta_i^s\| &\leq& \left( \|\bar{w}_i\| + 2 \, |\bar{\beta}_i| \, \|\bar{s}_{i-1}\| \right) \epsilon , \notag \\
	\|\delta_i^w\| &\leq& \left( \|\bar{w}_i\| + 2 \, |\bar{\alpha}_i| \, \|\bar{z}_i\| \right) \epsilon  , \notag \\
	\|\delta_i^q\| &\leq& \left( (\tilde{\mu}\sqrt{n}+1) \, \|M^{-1}\| \, \|\bar{w}_i\| + 2 \, |\bar{\beta}_i| \, \|\bar{q}_{i-1}\| \right) \epsilon \notag \\
	\|\delta_i^z\| &\leq& \left( (\mu\sqrt{n}+\tilde{\mu}\sqrt{n}+1)  \, \|A\| \, \|M^{-1}\| \, \|\bar{w}_i\| + 2 \, |\bar{\beta}_i| \, \|\bar{z}_{i-1}\| \right) \epsilon . 
	\label{eq:errorterms}
\end{eqnarray}
Here it is assumed that in the three-term recurrences for $\bar{r}_{i+1}$ and $\bar{u}_{i+1}$ in Eq.\,\eqref{eq:recur}, the second and third term are summed up first. The scalar $\mu$ denotes the maximum number of non-zeros in any row of $A$, and $\tilde{\mu}$ is the row-wise max.~nnz for $M^{-1}$.

In finite precision, the gap between the true (explicitly computed) residual $b - A \bar{x}_i$ and the recursive residual $\bar{r}_i$ is denoted 
\begin{equation}
	f_i := \left(b-A\bar{x}_i\right) - \bar{r}_i.
\end{equation}
For $i = 0$, the residual $\bar{r}_0$ is computed explicitly in Alg.\,\ref{algo::ppipe-cg-shift}, and the gap $f_0$ is the roundoff from computing $\bar{r}_0$ from $A$, $\bar{x}_0$ and $b$, i.e., $f_0 = b-A\bar{x}_0 - \text{fl}(b-A\bar{x}_0)$. The norm of this initial gap is bounded by
	$\|f_0\| \leq \left( (\mu\sqrt{n}+1) \, \|A\| \, \|\bar{x}_0\| + \|b\| \right) \epsilon.$
In iteration $i$ we obtain the following formula for the gap:
\begin{eqnarray} \label{eq:expression_f}
	f_{i+1} &=& (b-A\bar{x}_{i+1}) - \bar{r}_{i+1} \notag \\
					&=& b-A(\bar{x}_i + \bar{\alpha}_i \bar{p}_i + \delta_i^x) - (\bar{r}_i - \bar{\alpha}_i \bar{s}_i - \bar{\alpha}_i \sigma \bar{t}_i + \delta_i^r) \notag \\
					&=& f_i - \bar{\alpha}_i g_i - A\delta_i^x - \delta_i^r, \label{eq:recur_f_pipecg}
\end{eqnarray}
where $g_i$ is the gap between the true and recursive auxiliary vector $s_i$, i.e.,
\begin{equation}
	g_i := \left( A \bar{p}_i - \sigma \bar{t}_i \right) - \bar{s}_i.
\end{equation}
For $i = 0$ it holds that
	$\|g_0\| \leq \left( (\mu\sqrt{n}+1) \, \|A\| \, \|\bar{p}_0\| + 2 \, |\sigma| \, \| \bar{t}_0 \| \right) \epsilon.$
The residual gap in iteration $i$ is coupled to the error $g_i$, which can be written as
\begin{eqnarray}\label{eq:expression_g}
	g_i &=& \left( A \bar{p}_i - \sigma \bar{t}_i \right) - \bar{s}_i \notag \\
			&=& A (\bar{u}_i + \bar{\beta}_i \bar{p}_{i-1} +\delta_i^p) - \sigma (\bar{r}_i + \bar{\beta}_i \bar{t}_{i-1} + \delta^{t}_i) - (\bar{w}_i + \bar{\beta}_i \bar{s}_{i-1} + \delta_i^s) \notag \\
			&=& h_i + \bar{\beta}_i g_{i-1} + A \delta_i^p - \sigma \delta^{t}_i - \delta_i^s, \label{eq:recur_g_pipecg}
\end{eqnarray}
The auxiliary variable $\bar{w}_i$ is also computed recursively, resulting in a gap 
\begin{equation}
	h_{i+1} := \left( A \bar{u}_{i+1} - \sigma \bar{r}_{i+1} \right) - \bar{w}_{i+1}.
\end{equation}
For $i = 0$ one obtains the bound
	$\|h_0\| \leq \left( (\mu\sqrt{n}+1) \, \|A\| \, \|\bar{u}_0\| + 2 \, |\sigma| \, \|\bar{r}_0\| \right) \epsilon,$
whereas in iteration $i$ it holds that
\begin{eqnarray} \label{eq:expression_h}
	h_{i+1} &=& \left( A \bar{u}_{i+1} - \sigma \bar{r}_{i+1} \right) - \bar{w}_{i+1} \notag \\
					&=& A (\bar{u}_i - \bar{\alpha}_i \bar{q}_i - \bar{\alpha}_i \sigma \bar{p}_i +\delta_i^u) - \sigma (\bar{r}_i - \bar{\alpha}_i \bar{s}_i - \bar{\alpha}_i \sigma \bar{t}_i + \delta^r_i) - (\bar{w}_i - \bar{\alpha}_i \bar{z}_i + \delta_i^w) \notag \\
					&=& h_i - \bar{\alpha}_i j_i - \bar{\alpha}_i \sigma g_i + A \delta_i^u - \sigma \delta^r_i - \delta_i^w, \label{eq:recur_h_pipecg}
\end{eqnarray}
where the gap $j_i$ between the true and recursive variable $\bar{z}_i$ is 
\begin{equation}
	j_i = A \bar{q}_i - \bar{z}_i.
\end{equation}
For $i=0$ we can bound the norm of $j_i$ as the roundoff, i.e.,
	$\|j_0\| \leq \mu\sqrt{n} \, \|A\| \, \|\bar{q}_0\| \, \epsilon,$
while in iteration $i$ it holds that
\begin{eqnarray} \label{eq:expression_j}
	j_i &=& A \bar{q}_i - \bar{z}_i \notag \\
			&=& A (\bar{m}_i + \bar{\beta}_i \bar{q}_{i-1} + \delta_i^q) - (\bar{n}_i + \bar{\beta}_i \bar{z}_{i-1} + \delta_i^z) \notag \\
			&=& \bar{\beta}_i j_{i-1} + A \delta_i^q - \delta_i^z. \label{eq:recur_j_pipecg}
\end{eqnarray}
The final equation in Eq.\,\eqref{eq:recur_j_pipecg} holds since $\bar{n}_i = A\bar{m}_i$ is computed explicitly in Alg.\,\ref{algo::ppipe-cg-shift}.

\subsection{Local rounding error propagation matrices}

From Eq.\,\eqref{eq:expression_f}, Eq.\,\eqref{eq:expression_g}, Eq.\,\eqref{eq:expression_h} and Eq.\,\eqref{eq:expression_j} it follows that the residual gap for shifted pipelined CG is given by the following system of coupled equations:
\begin{equation} \label{eq:pipecg_system}
\begin{bmatrix}
 f_{i+1}  \\  g_i \\ h_{i+1} \\ j_i
\end{bmatrix} = 
\begin{bmatrix}
    1 & -\bar{\alpha}_i \bar{\beta}_i & -\bar{\alpha}_i & 0 \\ 0 & \bar{\beta}_i & 1 & 0 \\ 0 & - \bar{\alpha}_i \bar{\beta}_i \sigma & 1 - \bar{\alpha}_i \sigma & -\bar{\alpha}_i \bar{\beta}_i \\ 0 & 0 & 0 & \bar{\beta}_i
\end{bmatrix}
\begin{bmatrix}
f_i \\ g_{i-1} \\ h_i \\ j_{i-1}
\end{bmatrix} +
\begin{bmatrix}
	\epsilon^f_i\\
	\epsilon^g_i\\
	\epsilon^h_i\\
	\epsilon^j_i
\end{bmatrix},
\end{equation}
where the local rounding errors are
\begin{equation} 
\begin{bmatrix}
	\epsilon^f_i\\
	\epsilon^g_i\\
	\epsilon^h_i\\
	\epsilon^j_i
\end{bmatrix}
=
\begin{bmatrix}
  - A\delta_i^x - \delta_i^r - \bar{\alpha}_i \left( A \delta_i^p - \sigma \delta^{t}_i - \delta_i^s \right) \\ 
	A \delta_i^p - \sigma \delta^{t}_i - \delta_i^s \\
	A\delta_i^u - \sigma \delta^r_i- \delta_i^w - \bar{\alpha}_i \left( A \delta_i^q - \delta_i^z \right) - \bar{\alpha}_i \sigma \left( A \delta^p_i - \sigma \delta^{t}_i - \delta^s_i \right) \\ 
	A \delta_i^q - \delta_i^z 
\end{bmatrix}.
\end{equation}
This system can be written for short as
\begin{equation} \label{eq:system_short}
	\pi_{i+1} = P_i(\sigma) \pi_i + \epsilon_i,
\end{equation}
where $\pi_{i+1} = [f_{i+1}, g_i, h_{i+1}, j_i]^T$, the local error propagation matrix $P_i(\sigma)$ is
\begin{equation} \label{eq:propmat_shifted}
P_i(\sigma) =
\begin{bmatrix}
    1 & -\bar{\alpha}_i \bar{\beta}_i & -\bar{\alpha}_i & 0 \\ 0 & \bar{\beta}_i & 1 & 0 \\ 0 & - \bar{\alpha}_i \bar{\beta}_i \sigma & 1 - \bar{\alpha}_i \sigma & -\bar{\alpha}_i \bar{\beta}_i \\ 0 & 0 & 0 & \bar{\beta}_i
\end{bmatrix},
\end{equation}
and the local rounding errors in each iteration are $\epsilon_i = [\epsilon^f_i,\epsilon^g_i,\epsilon^h_i,\epsilon^j_i ]^T$.
Note that when the shift is zero, 
the local error propagation matrix $P_i(\sigma)$ reduces to
the local error propagation matrix for pipelined CG, Alg.\,\ref{algo::ppipe-cg}, which was derived in \cite{cools2016analysis}.

Via recursive substitution of Eq.\,\eqref{eq:system_short}, we find that the gaps on the variables $\bar{r}_{i+1}$, $\bar{s}_i$, $\bar{w}_{i+1}$ and $\bar{z}_i$ after $i$ iterations are given by
\begin{equation} \label{eq:system_unwinded}
	\pi_{i+1} = \left( \prod_{k=1}^i P_k(\sigma) \right) \pi_1 + \sum_{j=1}^{i} \left( \prod_{k=j+1}^i P_k(\sigma) \right) \epsilon_j,
\end{equation}
where 
the product notation in Eq.\,\eqref{eq:system_unwinded} should be interpreted as follows:
\begin{align}
	\mathbf{P}_{j,i}(\sigma) &:= \left(\prod_{k=j}^i P_k(\sigma)\right) := P_i(\sigma) \cdot P_{i-1}(\sigma) \cdot \ldots \cdot P_j(\sigma), &\text{for~} j \leq i, \notag \\
	\mathbf{P}_{j,i}(\sigma) &:= \left(\prod_{k=j}^i P_k(\sigma)\right) := I, &\text{for~} j > i. \label{eq:matrix_product}
\end{align} 
With the short-hand notation $\mathbf{P}_{j,i}(\sigma)$, introduced in Eq.\,\eqref{eq:matrix_product}, for the product of the local error propagation matrices, Eq.\,\eqref{eq:system_unwinded} becomes
\begin{equation} \label{eq:system_unwinded2}
	\pi_{i+1} = \mathbf{P}_{1,i}(\sigma) \, \pi_1 + \sum_{j=1}^i \mathbf{P}_{j+1,i}(\sigma) \, \epsilon_j.
\end{equation}
The propagation matrices $\mathbf{P}_{1,i}(\sigma), \mathbf{P}_{2,i}(\sigma), \ldots, \mathbf{P}_{i,i}(\sigma)$ thus fully characterize the gaps on the variables $\bar{r}_{i+1}$, $\bar{s}_i$, $\bar{w}_{i+1}$ and $\bar{z}_i$ in iteration $i+1$. 

In the case of classic CG, the local error propagation matrix $P_i(\sigma)$ is reduced to the scalar $1$, see standard rounding error analysis of CG in for example \cite{cools2016analysis,greenbaum1997estimating,gutknecht2000accuracy}. Indeed, in Alg.\,\ref{algo::pcg}, the auxiliary variables $s_i$, $w_i$ and $z_i$ are not computed (recursively), and it holds that
\begin{equation}\label{eq:gap_f0_cg}
	f_{i+1} = f_0  - \sum_{j=0}^i \left(A\delta_j^x + \delta_j^r\right) = f_0  - \sum_{j=0}^i \epsilon_i^f,
\end{equation}
and hence local rounding errors on the residual are thus merely accumulated, and no propagation of local errors occurs. 

For shifted pipelined CG Eq.\,\eqref{eq:system_unwinded2} dictates that, in a given iteration $i$, the magnitude of the entries of each individual matrix $\mathbf{P}_{j,i}(\sigma)$ with $j \leq i$ indicate whether the corresponding local rounding error is amplified.  When the modulus of a matrix element is of $\mathcal{O}(1)$, the corresponding local rounding error is merely accumulated, and the solution accuracy is expected to be comparable to classic CG. However, when an entry is significantly (several orders of magnitude) larger than one, the local rounding error component is propagated, and may have a detrimental impact on the residual gap. Attainable accuracy may then be reduced significantly compared to classic CG. This phenomenon lies at the heart of the loss of attainable accuracy in the pipelined CG method, Alg.\,\ref{algo::ppipe-cg}, as described in \cite{cools2016analysis}. 
We characterize the propagation matrix $\mathbf{P}_{j,i}(\sigma)$ by its 2-norm, and define the following function:
\begin{equation} \label{eq:functionpsi}
	\psi_i(\sigma) =  \max_{1\leq j \leq i} \|\mathbf{P}_{j,i}(\sigma)\|_2.
\end{equation}
The function $\psi_i(\sigma)$ has to be minimized in function of the shift $\sigma$ in order to determine the optimal shift for a given problem in a given iteration $i$. By considering $i$ large enough, preferably beyond the stagnation point for classic CG, a shift can be determined that allows to achieve a final p-CG-sh solution of comparable accuracy to classic CG. This is illustrated by numerical experiments in Section \ref{sec:numerical}.

\newpage

\section{Shifted pipelined Conjugate Gradients with variable shift}

\subsection{Derivation of the algorithm}

In this section we investigate whether it is possible to define a pipelined CG algorithm with variable, iteration-dependent shift parameter $\sigma_i$. This would allow for more flexibility in the shift choice to minimize the function $\psi_i(\sigma)$ defined by Eq.\,\eqref{eq:functionpsi}. 
We show that it is possible to extend the shifted algorithm, Alg.\,\ref{algo::ppipe-cg-shift}, to a version with a variable shift, by adapting the recurrences for the auxiliary variables $s_i$, $w_{i+1}$, $q_i$ and $z_i$.

We again start from the recurrences 
\begin{equation}
	t_i = r_i + \beta_i t_{i-1} , \quad \text{and} \quad p_i = u_i + \beta_i p_{i-1}.
\end{equation}
The auxiliary variable $s_i$ is now defined as $s_i := \left(A M^{-1} - \sigma_i I\right) t_i = A p_i - \sigma_i t_i$. 
From this definition it follows that
\begin{eqnarray}
	s_i = A p_i - \sigma_i t_i = A p_i - \sigma_{i-1} t_i - (\sigma_{i}-\sigma_{i-1}) t_i
\end{eqnarray}
After substituting the recurrence for $p_i$ and $t_i$ in the right-hand side above, we obtain
\begin{eqnarray}
s_i &=& A \left(u_i + \beta_i p_{i-1}\right) - \sigma_{i-1} \left(r_i + \beta_i t_{i-1}\right) - (\sigma_{i}-\sigma_{i-1}) t_i \notag \\
&=& w_i + \beta_i s_{i-1} - (\sigma_{i}-\sigma_{i-1}) t_i
\end{eqnarray}
where $w_i := \left(A M^{-1} - \sigma_{i-1} I \right) r_i = Au_i - \sigma_{i-1} r_i$. The recurrence for $s_i$ features a correction term $(\sigma_{i}-\sigma_{i-1}) t_i$ to account for the difference in the shift between successive iterations. Since $t_i$ is required to compute this term, $s_i$ can only be computed after the recursion for $t_i$ in the algorithm. For the residual we derive the recurrence
\begin{eqnarray}
	r_{i+1} &=& r_i - \alpha_i A p_i \notag \\
					&=& r_i - \alpha_i \left(s_i + \sigma_i t_i\right) \notag \\
					&=& r_i - \alpha_i s_i - \alpha_i \sigma_i t_i,
\end{eqnarray}
which is the exact analogue of Eq.\,\eqref{eq:r}, and for the preconditioned residual we have
\begin{eqnarray}
	u_{i+1} &=& u_i - \alpha_i M^{-1} s_i - \alpha_i \sigma_i M^{-1}t_i \notag \\
					&=& u_i - \alpha_i q_i - \alpha_i \sigma_i p_i,
\end{eqnarray}
similar to Eq.\,\eqref{eq:u}, where $q_i := M^{-1} s_i$ and $p_i := M^{-1} t_i$.
Multiplication of the recurrence for $r_{i+1}$ by $\left(A M^{-1} - \sigma_i I \right)$ on both sides yields
\begin{eqnarray}
\left(A M^{-1} - \sigma_i I \right) r_{i+1} &=& \left(A M^{-1} - \sigma_i I \right) \left( r_i - \alpha_i s_i - \alpha_i \sigma_i t_i \right) \notag \\
w_{i+1}	&=& \left(A M^{-1} - \sigma_{i-1} I \right) r_i - (\sigma_i-\sigma_{i-1}) r_i \notag \\
				& & ~ - \alpha_i \left(A M^{-1} - \sigma_i I \right) s_i - \alpha_i \sigma_i \left(A M^{-1} - \sigma_i I \right) t_i \notag \\
				&=& w_i - \alpha_i z_i - (\sigma_i - \sigma_{i-1}) r_i,
\end{eqnarray}
where $z_i = A q_i = A M^{-1} s_i$. Thus, also in the recurrence for $w_{i+1}$ a correction has to be made based on the difference between consecutive shifts. Due to the term involving $r_i$, the recurrence for $w_{i+1}$ has to be computed before the recurrence for $r_{i+1}$, which overwrites the residual.
By multiplying the recurrence for $s_i$ by $M^{-1}$ we obtain 
\begin{eqnarray}
	M^{-1}s_i &=& M^{-1} w_i + \beta_i M^{-1} s_{i-1} - (\sigma_i-\sigma_{i-1}) M^{-1} t_i \notag \\
				q_i &=& m_i + \beta_i q_{i-1} - (\sigma_i-\sigma_{i-1}) p_i,
\end{eqnarray}
with $m_i = M^{-1} w_i$, and by multiplying this recursion for $q_i$ by $A$ we get
\begin{eqnarray}
	A q_i &=& A m_i + \beta_i A q_{i-1} - (\sigma_i-\sigma_{i-1}) A p_i \notag \\
		z_i &=& n_i + \beta_i z_{i-1}- (\sigma_i-\sigma_{i-1}) (s_i+\sigma_it_i).
\end{eqnarray}
We stress that $q_i$ cannot be computed until after the recurrence for $p_i$, since the latter variable is used in the recurrence for $q_i$. Likewise, $z_i$ can only be computed after the recurrences for $s_i$ and $t_i$ have been computed. The resulting pipelined CG method with variable shift is shown in Alg.\,\ref{algo::ppipe-cg-var-shift}.

\begin{algorithm}[t]
  \caption{Variable shifted preconditioned pipelined CG}
  \label{algo::ppipe-cg-var-shift}
  \begin{algorithmic}[1]
  	\Procedure{p-cg-var-sh}{$A$, $M^{-1}$, $b$, $x_0$, $\sigma_{-1}$, $\sigma_0$, $\sigma_1$, \ldots}
    \State $r_0 := b - Ax_0$; $u_0:= M^{-1} r_0$; $w_0 := Au_0 - \sigma_{-1} r_0$
    \For{$i = 0,\dots$}
    \State $\gamma_i :=(r_i,u_i)$
    \State $\delta := (w_i+\sigma_{i-1} r_i,u_i)$
    \State $m_i := M^{-1} w_i$
    \State $n_i := A m_i$
    \If{$i>0$}
    \State $\beta_i := \gamma_i/\gamma_{i-1}$; $\alpha_i := (\delta/\gamma_i - \beta_i/\alpha_{i-1})^{-1}$
    \Else
    \State $\beta_i :=0$; $\alpha_i := \gamma_i/\delta$
    \EndIf
		\State $t_i := r_i + \beta_i t_{i-1}$
    \State $p_i := u_i + \beta_i p_{i-1}$
		\State $s_i := w_i + \beta_i s_{i-1} - (\sigma_i - \sigma_{i-1}) t_i$
    \State $q_i := m_i + \beta_i q_{i-1} - (\sigma_i - \sigma_{i-1}) p_i$
		\State $z_i := n_i + \beta_i z_{i-1} - (\sigma_i - \sigma_{i-1}) (s_i+\sigma_it_i)$
    \State $x_{i+1} := x_i + \alpha_i p_i$
		\State $w_{i+1} := w_i - \alpha_i z_i - (\sigma_i-\sigma_{i-1}) r_i$
    \State $r_{i+1} := r_i - \alpha_i s_i - \alpha_i \sigma_i t_i$
    \State $u_{i+1} := u_i - \alpha_i q_i - \alpha_i \sigma_i p_i$
    \EndFor
    \EndProcedure
  \end{algorithmic}
\end{algorithm}

\subsection{Numerical stability analysis for variable shifted pipelined CG}

We adopt the notation from Section \ref{sec:rounding}. In finite precision the recurrences for the auxiliary variables in the variable shifted pipelined CG algorithm, Alg.\,\ref{algo::ppipe-cg-var-shift} are
\begin{align}
	\bar{x}_{i+1} &= \bar{x}_i + \bar{\alpha}_i \bar{p}_i + \delta_i^x, &
	\bar{t}_i &= \bar{r}_i + \bar{\beta}_i \bar{t}_{i-1} + \delta_i^{t} , \notag \\
	\bar{r}_{i+1} &= \bar{r}_i - \bar{\alpha}_i \bar{s}_i -\bar{\alpha}_i \sigma_i \bar{t}_i+ \delta_i^r,&
	\bar{p}_i 		&= \bar{u}_i + \bar{\beta}_i \bar{p}_{i-1} + \delta_i^p , \notag \\
	\bar{u}_{i+1} &= \bar{u}_i - \bar{\alpha}_i \bar{q}_i - \bar{\alpha}_i \sigma_i \bar{p}_i + \delta_i^u , &
	\bar{s}_i 		&= \bar{w}_i + \bar{\beta}_i \bar{s}_{i-1} - \left(\sigma_i-\sigma_{i-1}\right) \bar{t}_i + \delta_i^s, \notag \\
	\bar{w}_{i+1} &= \bar{w}_i - \bar{\alpha}_i \bar{z}_i - \left(\sigma_i-\sigma_{i-1}\right) \bar{r}_i + \delta_i^w, &	
	\bar{q}_i     &= \bar{m}_i + \bar{\beta}_i \bar{q}_{i-1} - \left(\sigma_i-\sigma_{i-1}\right) \bar{p}_i + \delta_i^q, \notag \\
	& &
	\bar{z}_i     &= \bar{n}_i + \bar{\beta}_i \bar{z}_{i-1} - \left(\sigma_i-\sigma_{i-1}\right) \left(\bar{s}_i - \sigma_i \bar{t_i} \right) + \delta_i^z, \label{eq:recur2}
\end{align}
where the local rounding errors on each variable are bounded as follows
\begin{align} 
	\|\delta_i^x\| &\leq \left( \|\bar{x}_i\| + 2 \, |\bar{\alpha}_i| \, \|\bar{p}_i\| \right) \epsilon , \notag \\
	\|\delta_i^{t}\| &\leq \left( \|\bar{r}_i\| + 2 \, |\bar{\beta}_i| \, \|\bar{t}_{i-1}\| \right) \epsilon , \notag \\
	\|\delta_i^r\| &\leq \left( \|\bar{r}_i\| + 3 \, |\bar{\alpha}_i| \, \|\bar{s}_i\| + 4 \, |\bar{\alpha}_i| \, |\sigma_i| \, \|\bar{t}_i\| \right) \epsilon , \notag \\
	\|\delta_i^p\| &\leq \left( \|\bar{u}_i\| + 2 \, |\bar{\beta}_i| \, \|\bar{p}_{i-1}\| \right) \epsilon , \notag \\
	\|\delta_i^u\| &\leq \left( \|\bar{u}_i\| + 3 \, |\bar{\alpha}_i| \, \|\bar{q}_i\| + 4 \, |\bar{\alpha}_i| \, |\sigma_i| \, \|\bar{p}_i\| \right) \epsilon  , \notag \\
	\|\delta_i^s\| &\leq \left( \|\bar{w}_i\| + 3 \, |\bar{\beta}_i| \, \|\bar{s}_{i-1}\| + 4 \, |\sigma_i - \sigma_{i-1}| \, \|\bar{t}_i\| \right) \epsilon , \notag \\
	\|\delta_i^w\| &\leq \left( \|\bar{w}_i\| + 3 \, |\bar{\alpha}_i| \, \|\bar{z}_i\| + 4 \,|\sigma_i-\sigma_{i-1}| \, \|\bar{r}_i\| \right) \epsilon  , \notag \\
	\|\delta_i^q\| &\leq \left( (\tilde{\mu}\sqrt{n}+1) \, \|M^{-1}\| \, \|\bar{w}_i\| + 3 \, |\bar{\beta}_i| \, \|\bar{q}_{i-1}\| +  4 \, |\sigma_i - \sigma_{i-1} | \, \|\bar{p_i}\| \right) \epsilon \notag \\
	\|\delta_i^z\| &\leq \left( (\mu\sqrt{n}+\tilde{\mu}\sqrt{n}+1)  \, \|A\| \, \|M^{-1}\| \, \|\bar{w}_i\| + 3 \, |\bar{\beta}_i| \, \|\bar{z}_{i-1}\| + \ldots \right. \notag \\
	& ~~~~ \left. \ldots + 5 \, |\sigma_i-\sigma_{i-1}| \, \|\bar{s}_i\| + 6 \, |\sigma_i-\sigma_{i-1}| \, |\sigma_i| \, \|\bar{t_i}\| \right) \epsilon. \label{eq:errorterms2}
\end{align}
The gap between the true and recursive residual, $f_{i+1} = \left(b-A\bar{x}_i\right) - \bar{r}_i$, satisfies 
\begin{equation}
	f_{i+1} = f_i - \bar{\alpha}_i g_i - A\delta_i^x - \delta_i^r,
\end{equation}
which is identical to Eq.\,\eqref{eq:recur_f_pipecg}, with $g_i = \left( A \bar{p}_i - \sigma_i \bar{t}_i \right) - \bar{s}_i$. For the gap $g_i$ one obtains the following recurrence relation:
\begin{eqnarray}
	g_i &=& \left( A \bar{p}_i - \sigma_i \bar{t}_i \right) - \bar{s}_i \notag \\
			&=& A (\bar{u}_i + \bar{\beta}_i \bar{p}_{i-1} +\delta_i^p) - \sigma_i (\bar{r}_i + \bar{\beta}_i \bar{t}_{i-1} + \delta^{t}_i) \notag \\
			& & ~ - \left(\bar{w}_i + \bar{\beta}_i \bar{s}_{i-1} - (\sigma_i - \sigma_{i-1}) \bar{t}_i+ \delta_i^s\right) \notag \\
			&=& h_i + \bar{\beta}_i g_{i-1} + A \delta_i^p - \sigma_i \delta^{t}_i - \delta_i^s,
\end{eqnarray}
where $h_i = (A \bar{u}_i - \sigma_{i-1} \bar{r}_i) - \bar{w}_i$. This expression for $g_i$ is also identical to the case with fixed shift, see Eq.\,\eqref{eq:recur_g_pipecg}. Note that the bound on the local rounding error $\delta_i^s$ is different as before, since the recurrence for $\bar{s}_i$ has been modified.
The same remark can be made for $\delta_i^w$, $\delta_i^q$ and $\delta_i^z$.
Next, we derive the recurrence for $h_{i+1}$:
\begin{eqnarray}
	h_{i+1} &=& \left( A \bar{u}_{i+1} - \sigma_i \bar{r}_{i+1} \right) - \bar{w}_{i+1} \notag \\
					&=& A (\bar{u}_i - \bar{\alpha}_i \bar{q}_i - \bar{\alpha}_i \sigma_i \bar{p}_i +\delta_i^u) - \sigma_i (\bar{r}_i - \bar{\alpha}_i \bar{s}_i - \bar{\alpha}_i \sigma_i \bar{t}_i + \delta^r_i) \notag \\
					& & ~ - \left(\bar{w}_i - \bar{\alpha}_i \bar{z}_i - (\sigma_i - \sigma_{i-1}) \bar{r}_i + \delta_i^w \right) \notag \\
					&=& h_i - \bar{\alpha}_i j_i - \bar{\alpha}_i \sigma_i g_i + A \delta_i^u - \sigma_i \delta^r_i - \delta_i^w, 
\end{eqnarray}
with $j_i = A\bar{q}_i - \bar{z}_i$. Finally, for the gap $j_i$ we find 
\begin{eqnarray}
	j_i &=& A \bar{q}_i - \bar{z}_i \notag \\
			&=& A \left(\bar{m}_i + \bar{\beta}_i \bar{q}_{i-1} - (\sigma_i-\sigma_{i-1})\bar{p}_i + \delta_i^q \right) \notag \\
			& & - ~ \left(\bar{n}_i + \bar{\beta}_i \bar{z}_{i-1} - (\sigma_i-\sigma_{i-1})(\bar{s}_i+\sigma_i \bar{t}_i)+ \delta_i^z\right) \notag \\
			&=& \bar{\beta}_i j_{i-1} - (\sigma_i-\sigma_{i-1}) g_i + A \delta_i^q - \delta_i^z.
\end{eqnarray}
Consequently, we find that the error propagation matrix for variable shifted pipelined CG in iteration $i$ is 
\begin{equation} \label{eq:propmat_varshifted}
P_i(\sigma_{i-1},\sigma_i) = 
\begin{bmatrix}
    1 & -\bar{\alpha}_i \bar{\beta}_i & -\bar{\alpha}_i & 0 \\ 
		0 & \bar{\beta}_i & 1 & 0 \\ 
		0 & - \bar{\alpha}_i \bar{\beta}_i \sigma_{i-1} & 1 - \bar{\alpha}_i \sigma_{i-1} & -\bar{\alpha}_i \bar{\beta}_i \\ 
		0 & -\bar{\beta}_i(\sigma_i-\sigma_{i-1}) & -(\sigma_i-\sigma_{i-1}) & \bar{\beta}_i
\end{bmatrix}.
\end{equation}

\newpage

\section{Experimental results} \label{sec:numerical}

\subsection{Numerical accuracy results}

We present various examples from the Matrix Market\footnote{\url{http://math.nist.gov/MatrixMarket/}} library to illustrate the stabilizing properties of the shifted pipelined CG algorithm, Alg.\,\ref{algo::ppipe-cg-shift}. Details of four selected benchmark problems are given in Table \ref{tab:benchmarks_specs}. For all problems a significant loss of maximal attainable accuracy is observed when using the pipelined CG method, Alg.\,\ref{algo::ppipe-cg}, compared to classic CG, Alg.\,\ref{algo::pcg}. The right-hand side for all model problems is $b = A\hat{x}$ where $\hat{x}_j = 1/\sqrt{n}$, except for the \texttt{lapl200} benchmark where we use $b_j = 1/\sqrt{n}$. An all-zero initial guess $\bar{x}_0 = 0$ is used. Jacobi diagonal preconditioning (JAC) and Incomplete Cholesky Factorization (ICC) are included to reduce the number of Krylov iterations where required. For the preconditioner designated as $^*$ICC an compensated Incomplete Cholesky factorization is performed, where a real non-negative scalar $\eta = 0.5$ is used as a global diagonal shift in forming the Cholesky factor. The choice of the particular shift $\sigma^*$ is based on the numerical analysis in Section \ref{sec:numerical}. It is chosen in an \emph{a posteriori} fashion based on the history of the scalar coefficients $\alpha_i$ and $\beta_i$ computed by the pipelined CG method. 

\begin{figure}[t]
\begin{center}
\includegraphics[width=0.48\textwidth]{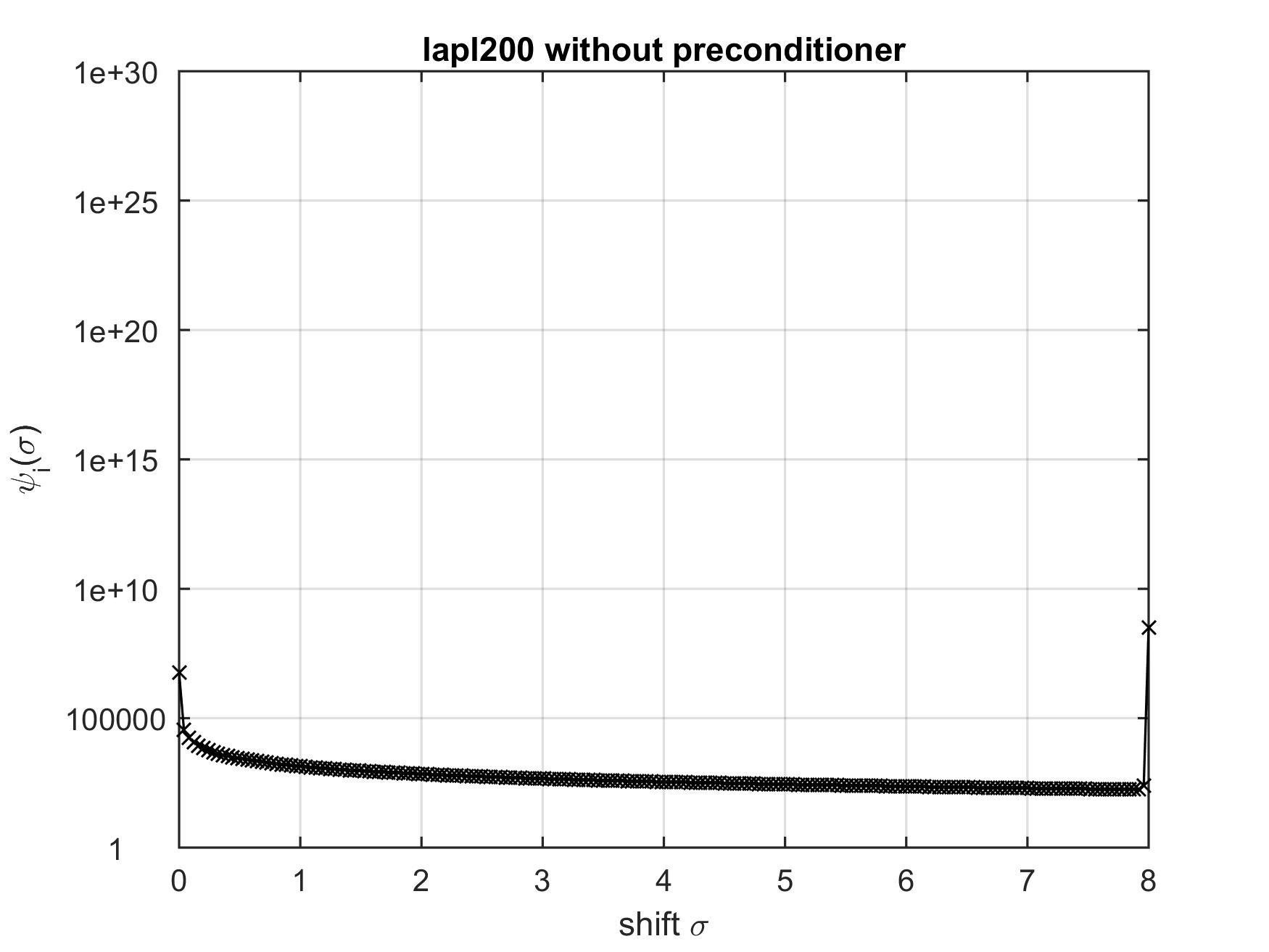}
\includegraphics[width=0.48\textwidth]{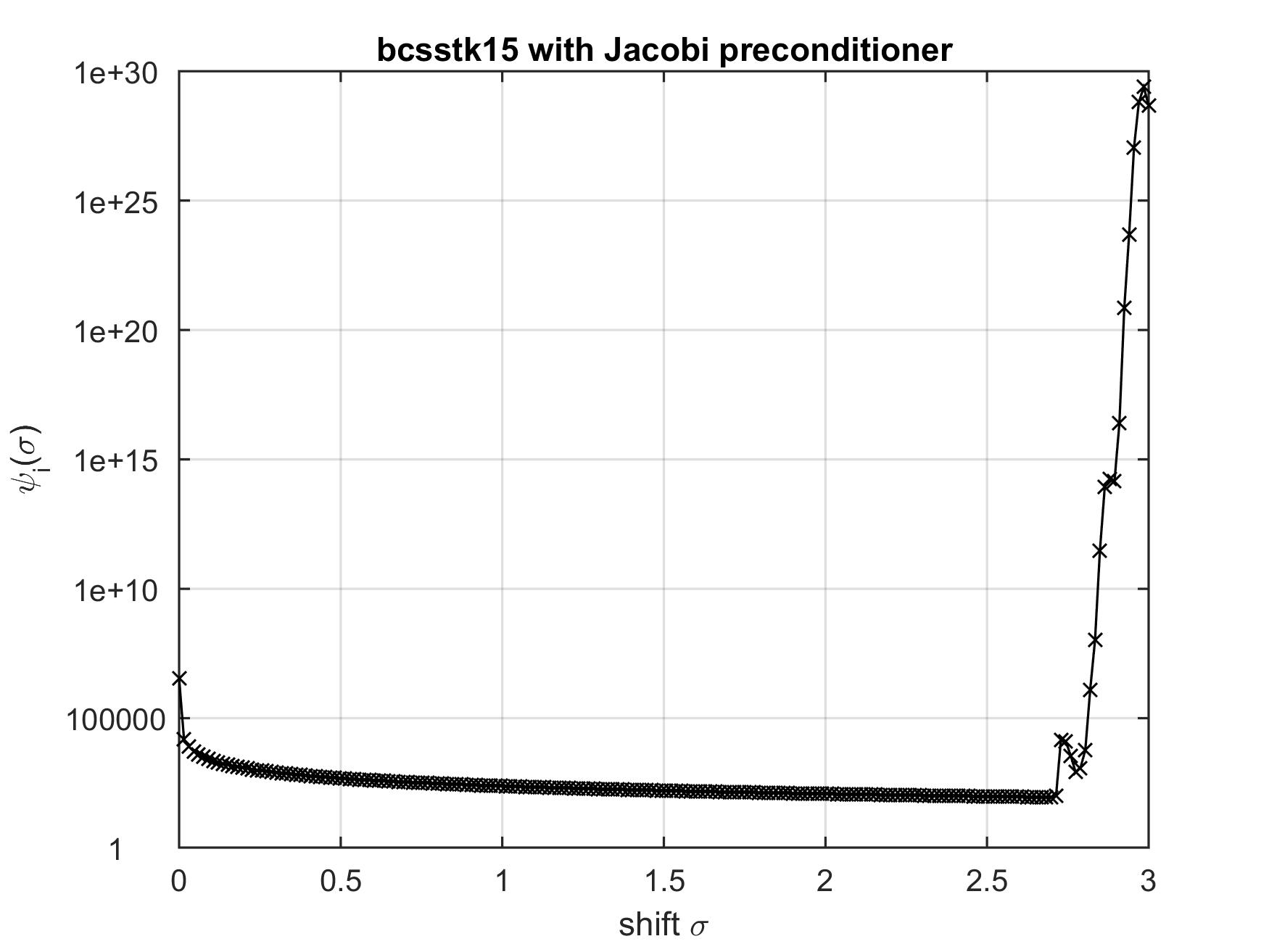} \\
\includegraphics[width=0.48\textwidth]{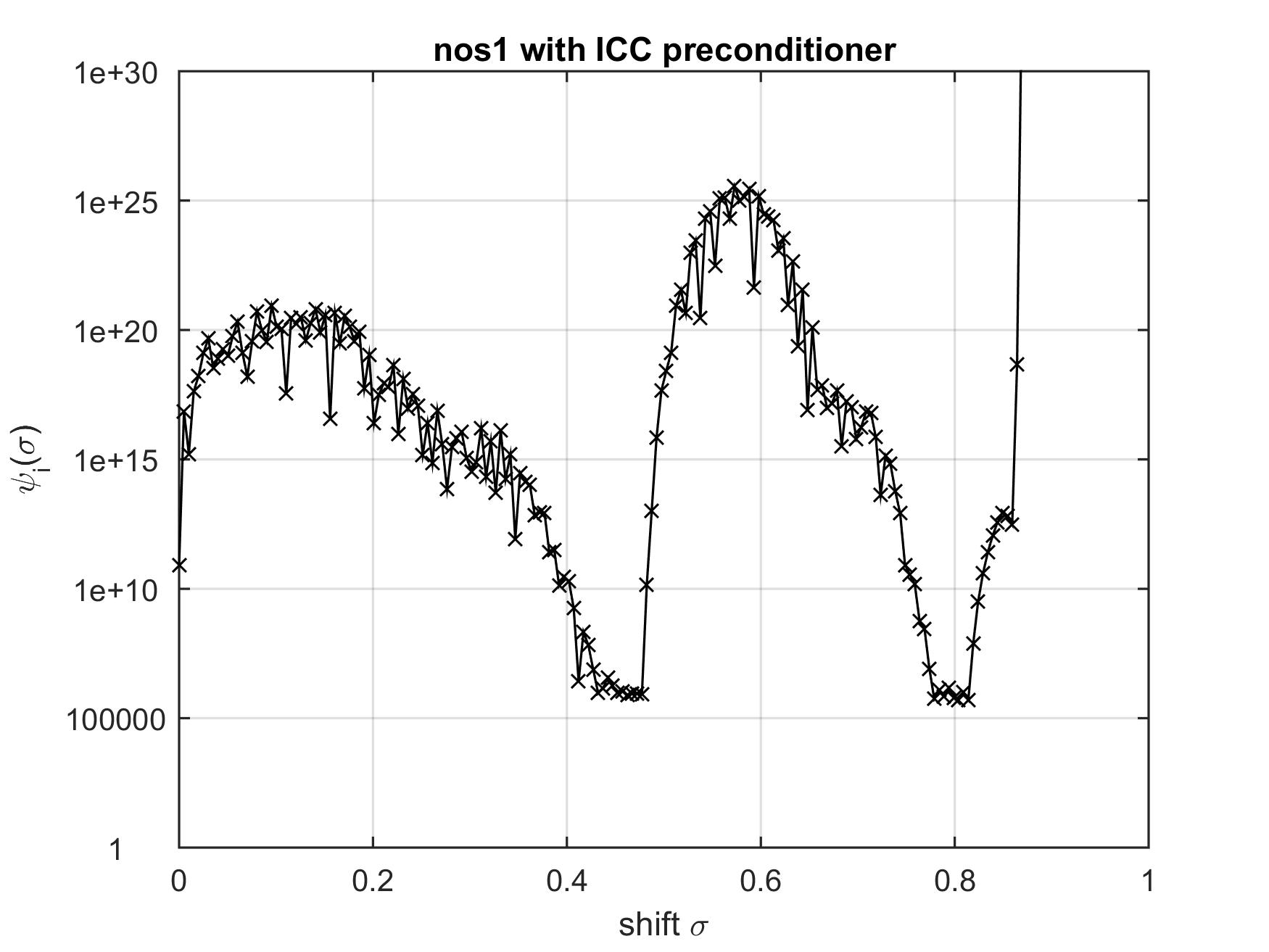}
\includegraphics[width=0.48\textwidth]{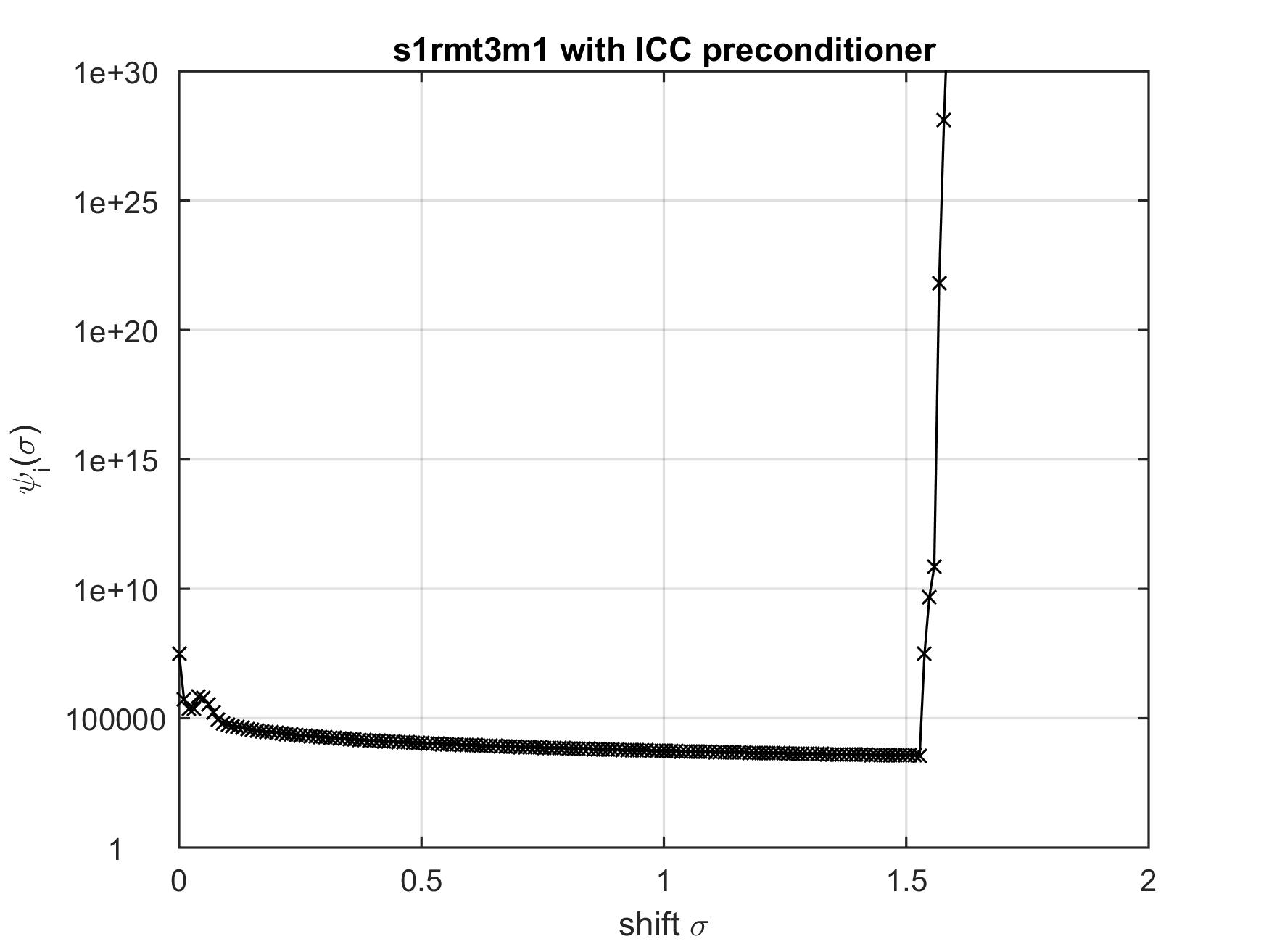}
\end{center}
\caption{Benchmark problems \texttt{lapl200} (top left), \texttt{bcsstk15} (top right), \texttt{nos1} (bottom left) and \texttt{s1rmt3m1} (bottom right). Evaluations of $\psi_i(\sigma)$ as a function of the shift choice $\sigma$ for a fixed number of iterations $i$ (see Table \ref{tab:benchmarks_specs}).}
\label{fig:shift}
\end{figure}

Figure \ref{fig:shift} shows discrete evaluations of the function $\psi_i(\sigma)$ as a function of the shift $\sigma$ for the maximum number of iterations $i$ listed in Table \ref{tab:benchmarks_specs}. The values of $\psi_i(\sigma)$ for standard, unshifted p-CG (i.e., with $\sigma = 0$) equal 1.78e+08 (\texttt{lapl200}), 4.34e+06 (\texttt{bcsstk15}), 8.30e+10 (\texttt{nos1}) and 1.12e+07 (\texttt{s1rmt3m1}) respectively. For the shifted p-CG method a shift $\sigma^*$ that stabilizes the method is chosen, see Table \ref{tab:benchmarks_specs}. The function $\psi_i$ is (close to) minimal for these shift choices and is several orders of magnitude smaller compared to the unshifted case, taking on the values 3.46e+02 (\texttt{lapl200}), 1.19e+02 (\texttt{bcsstk15}), 1.78e+08 (\texttt{nos1}) and 5.45e+02 (\texttt{s1rmt3m1}). Note that while choosing the shift too small inevitably results in the loss of accuracy displayed by the pipelined method, a (too) large value for the shift may destroy convergence entirely.

\begin{figure}[H]
\begin{center}
\includegraphics[width=0.48\textwidth]{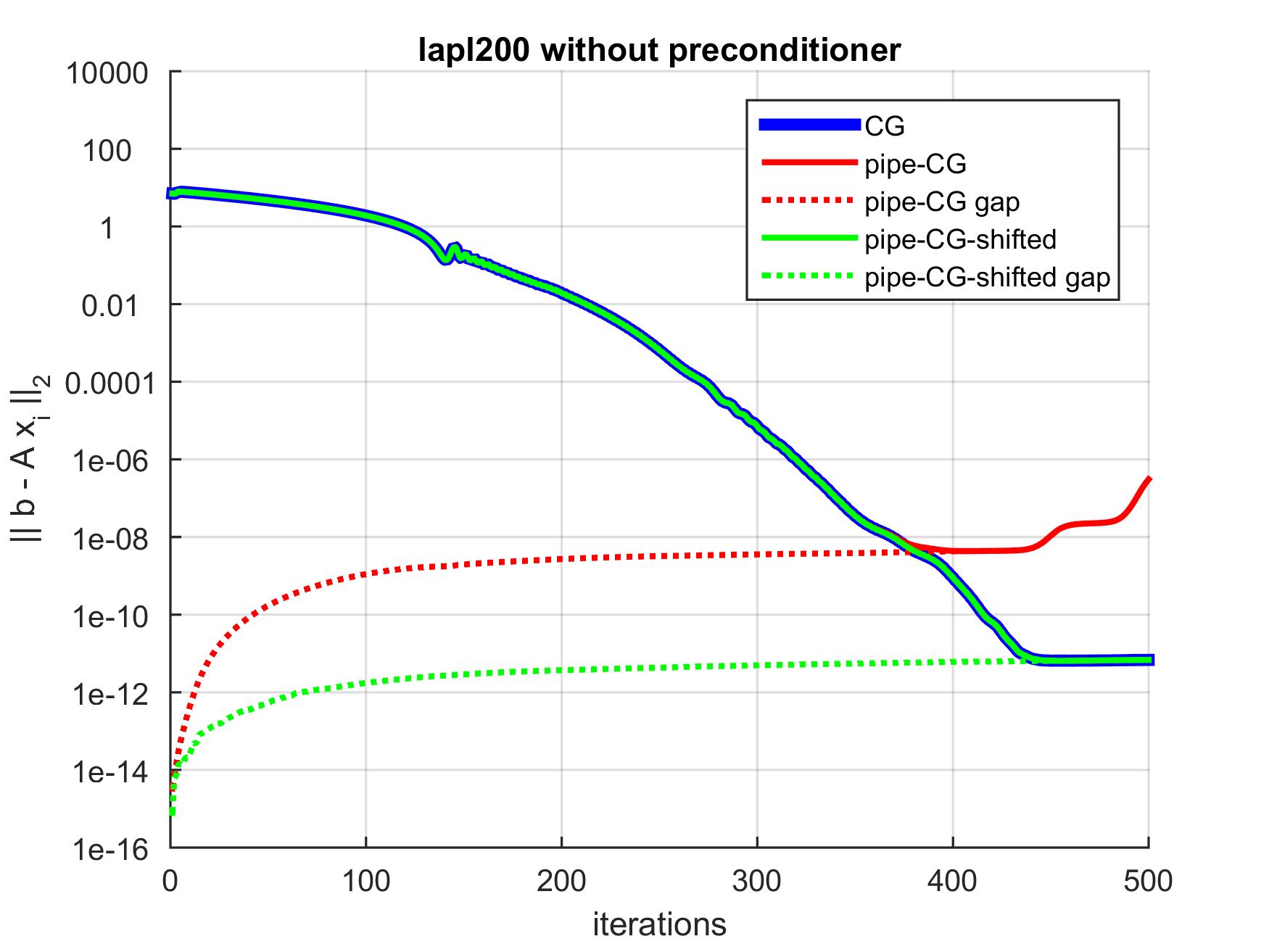}
\includegraphics[width=0.48\textwidth]{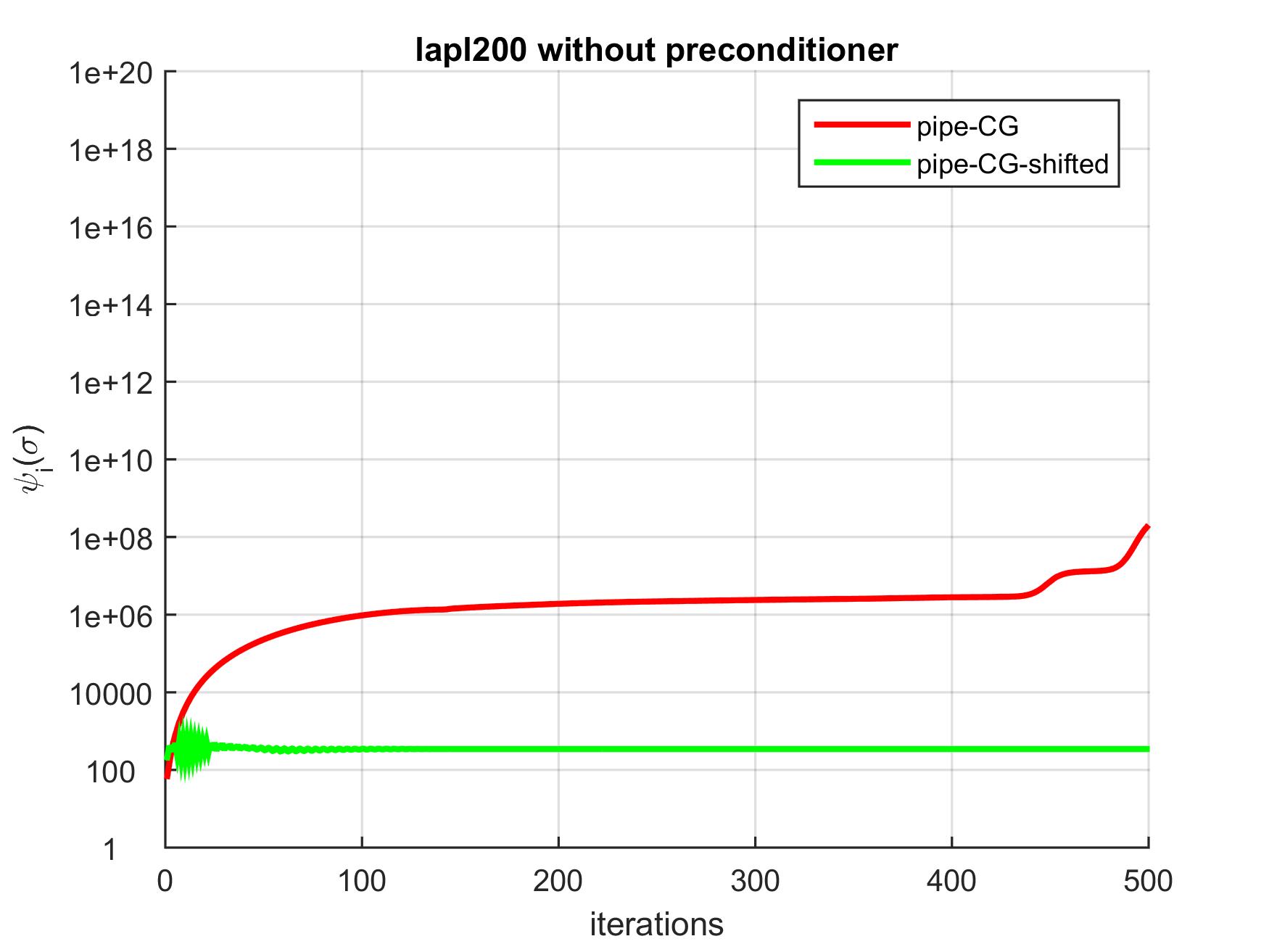}  \\
\includegraphics[width=0.48\textwidth]{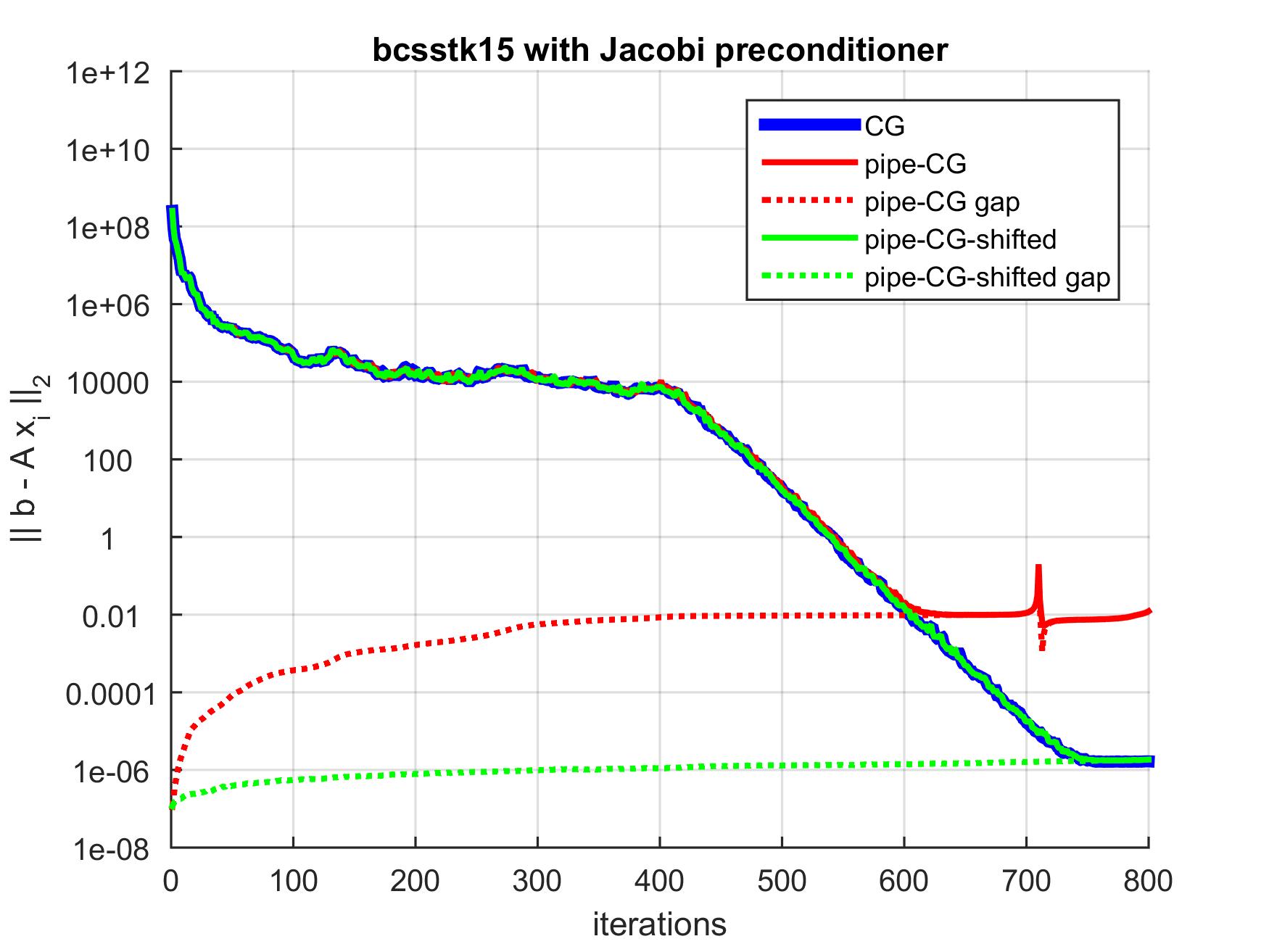}
\includegraphics[width=0.48\textwidth]{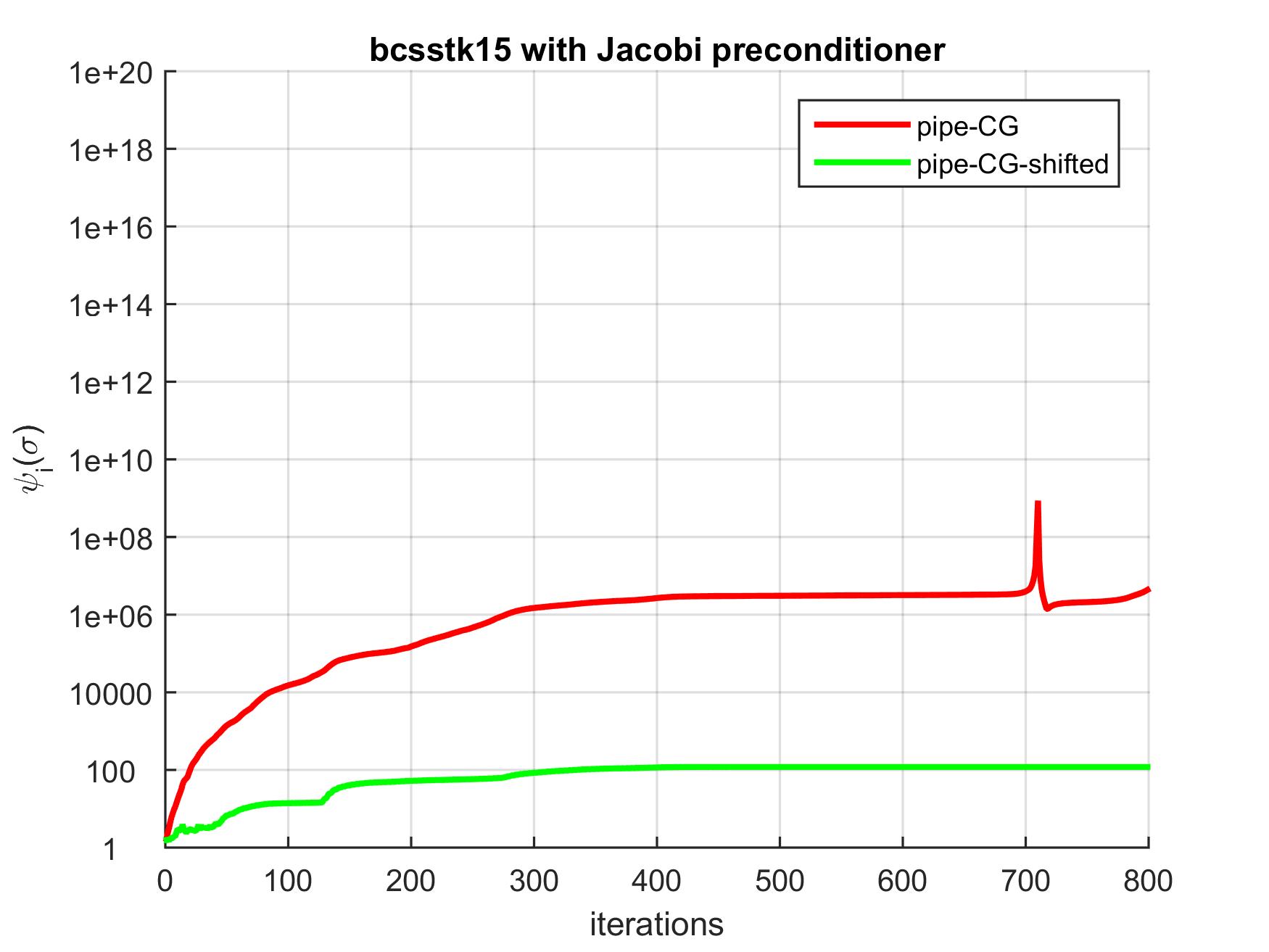}  \\ 
\includegraphics[width=0.48\textwidth]{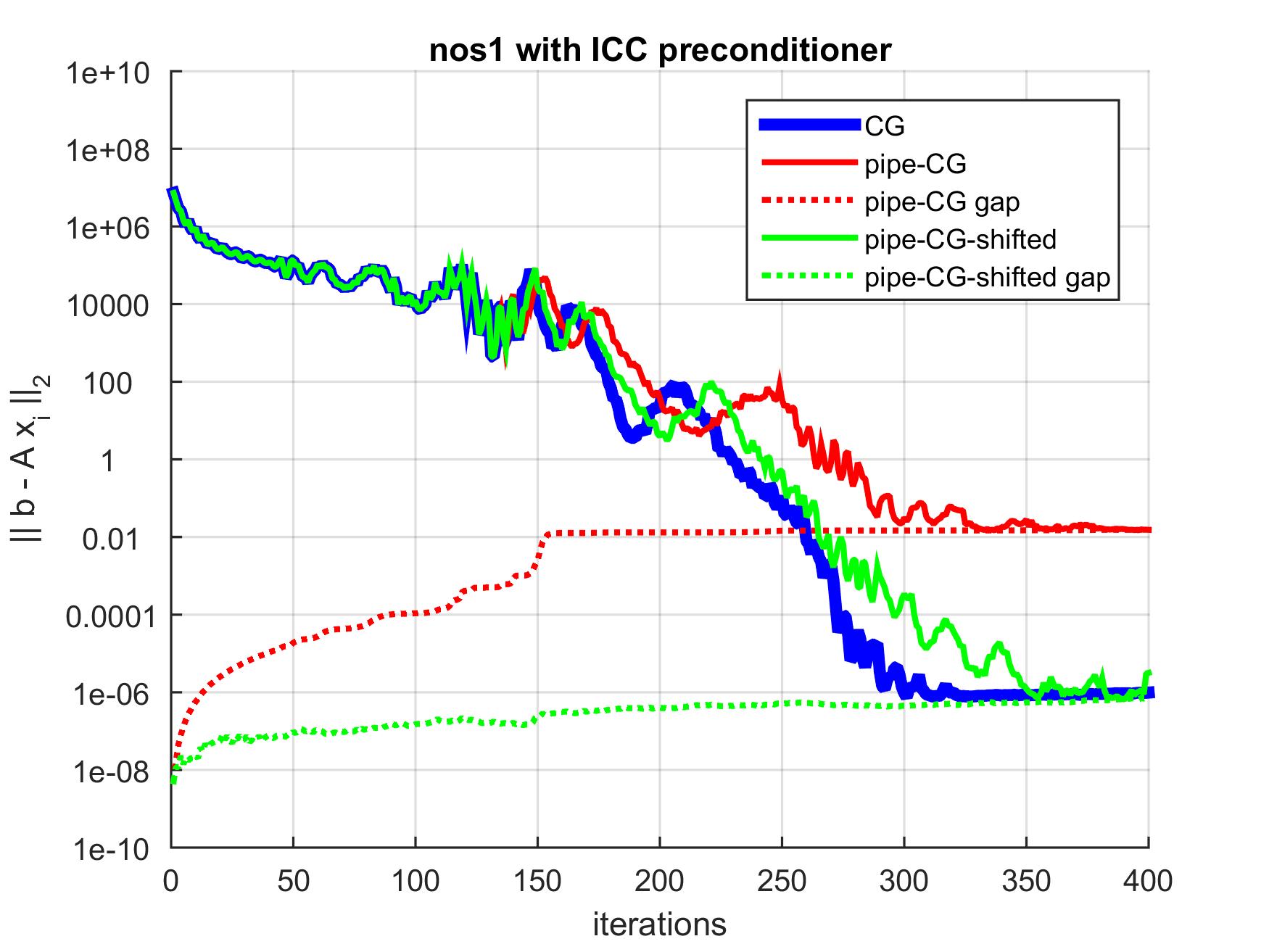}
\includegraphics[width=0.48\textwidth]{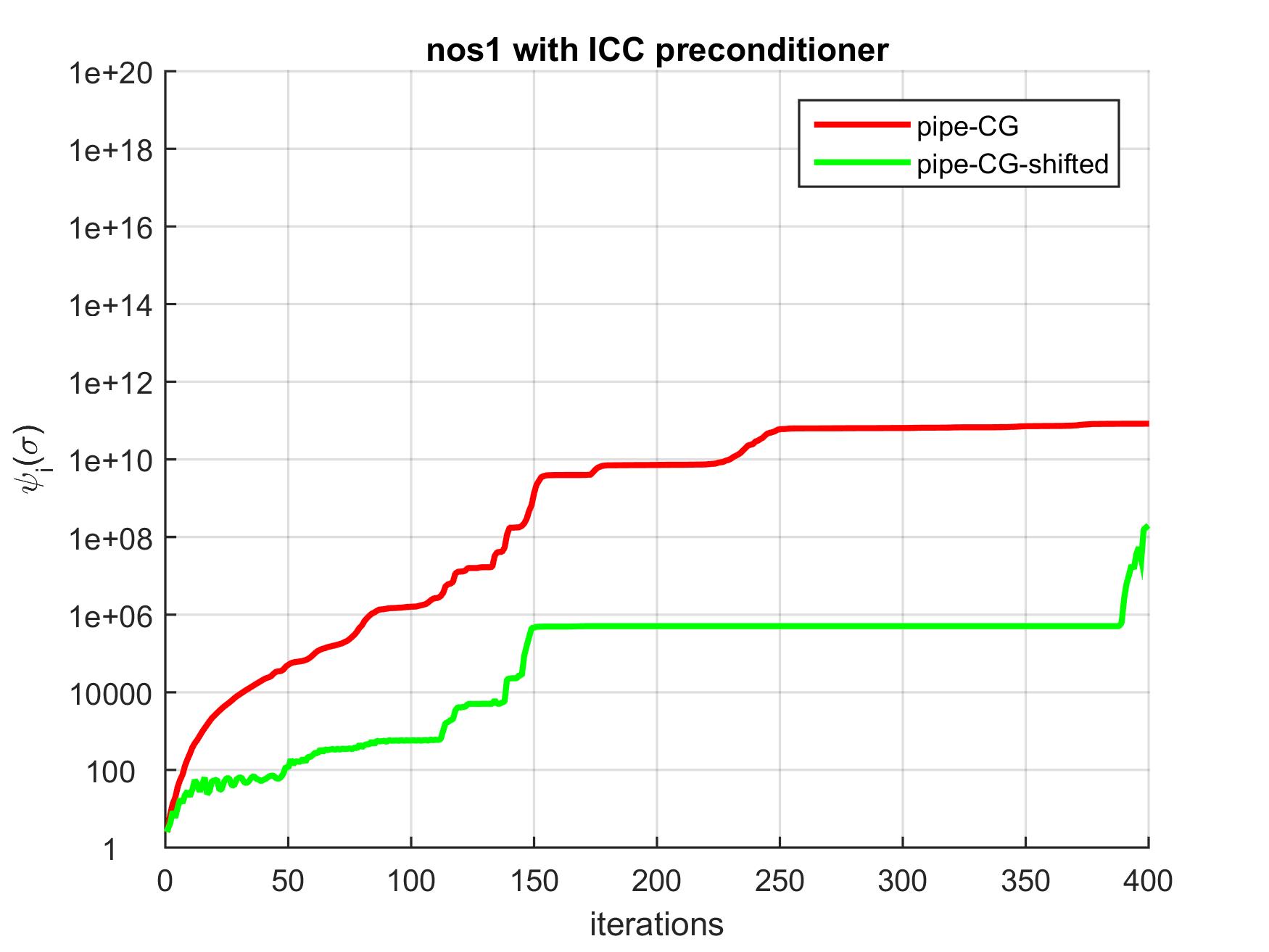}  \\ 
\includegraphics[width=0.48\textwidth]{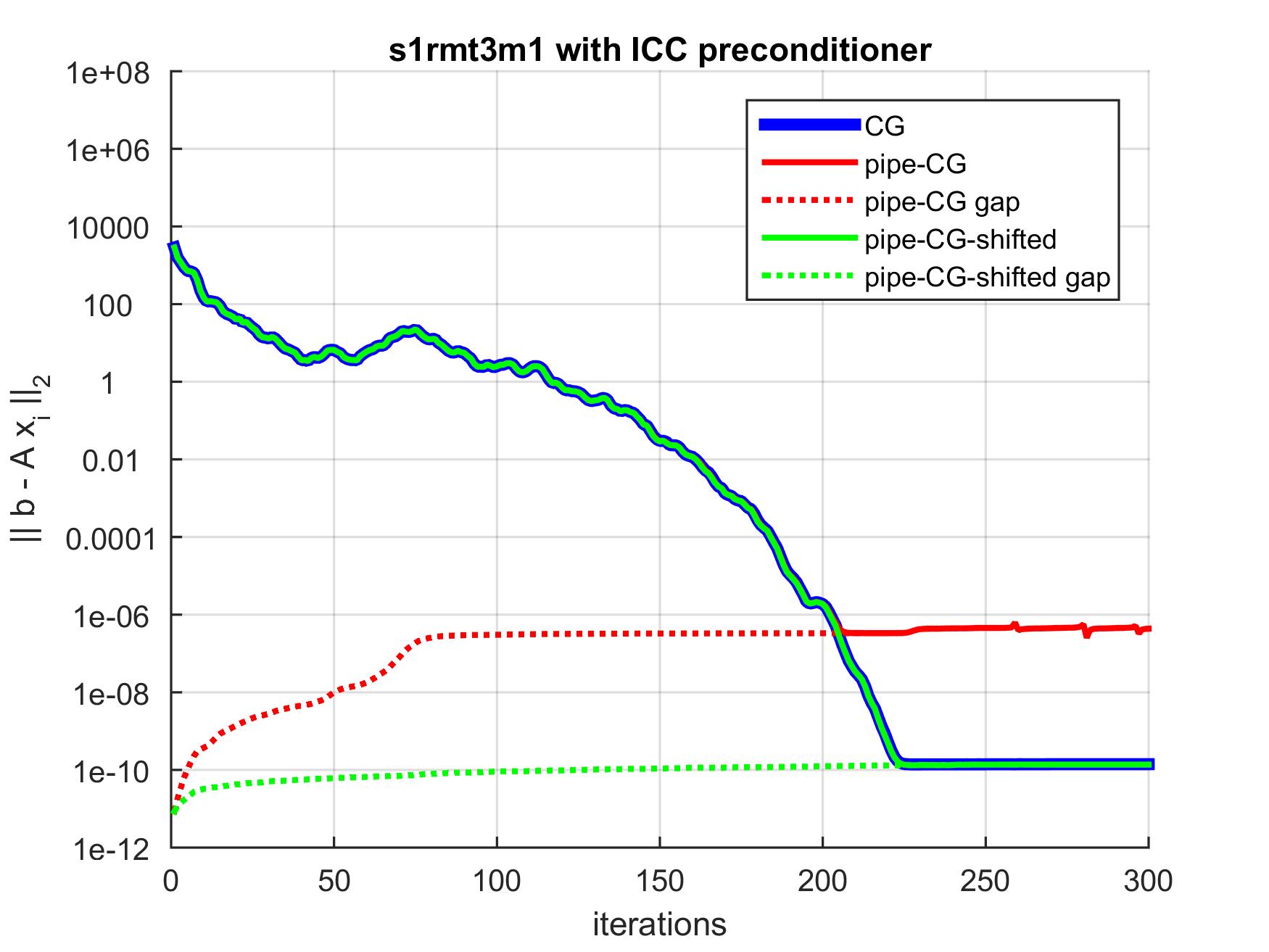}
\includegraphics[width=0.48\textwidth]{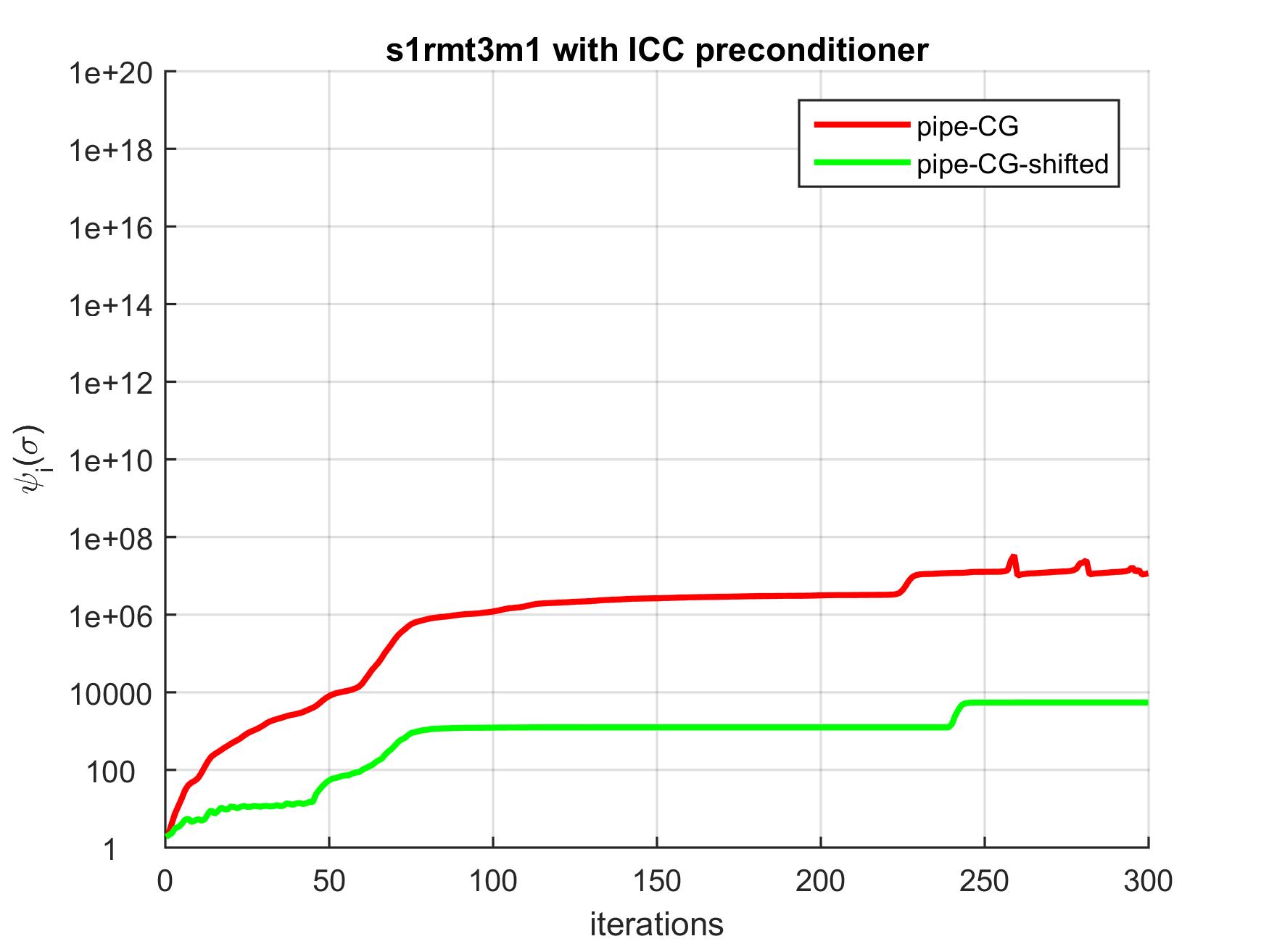}  \\  
\end{center}
\caption{Benchmark problems \texttt{lapl200}, \texttt{bcsstk15}, \texttt{nos1} and \texttt{s1rmt3m1}. Left: residual history $\|b-A\bar{x}_i\|_2$ (full line) and residual gap $\|(b-A\bar{x}_i)-\bar{r}_i\|_2$ (dotted line) as a function of iterations. Right: function $\psi_i(\sigma)$ as a function of iterations. The shift value is $\sigma = 0$ and $\sigma^*$ (see Table \ref{tab:benchmarks_specs}) for p-CG and shifted p-CG respectively.}
\label{fig:res1}
\end{figure}

Figures \ref{fig:res1} show the residual histories for the four benchmark problems (left panel) and the corresponding function $\psi_i(\sigma)$ as a function of iterations (right panel) evaluated in $\sigma = 0$ and $\sigma = \sigma^*$ (see Table \ref{tab:benchmarks_specs}) for the p-CG and shifted p-CG method respectively.  Note that the relative difference between the values $\psi_i(0)$ and $\psi_i(\sigma^*)$ (right panel) gives an estimate of the difference between the residual gaps (left panel, dotted lines), thus characterizing the improvement in attainable accuracy that is achieved by the shifted p-CG method. For all benchmark problems, the shifted p-CG method attains the accuracy of the classic CG method in a comparable number of iterations.

\begin{table}[t]
\centering
\scriptsize
\begin{tabular}{| l | r | r | r | r | r | r | r  r  r |}
\hline 
	Matrix 	 	& $n$    & Prec 		&	$\kappa(A)$	& $i$ & $\sigma^*$ &	$\|r_0\|_2$ 	& \multicolumn{3}{|c|}{$\|b-Ax_i\|_2$} \\
						&				 &					&					  	&				&				 &		& CG & p-CG & p-CG-$\sigma$ \\
\hline \hline
	\texttt{lapl200} 	& 40,000 & -  			& 2.4e+04  	& 500 & 4.00 &	1.0e+00 & 6.8e-12 & 3.1e-07 & 6.8e-12	\\
	\texttt{bcsstk15} 	&  3,948 & JAC  		& 8.0e+09 	& 800 & 2.00 &	4.3e+08	& 1.7e-06 & 1.2e-02 & 1.9e-06 \\
	\texttt{nos1} 			& 	 237 & $^*$ICC	& 2.5e+07 	& 400 & 0.82 &	5.7e+07	& 9.8e-07 & 1.5e-02 & 3.2e-06	\\
	\texttt{s1rmt3m1} 	&  5,489 & ICC  		& 2.5e+06 	& 300	& 1.00 &	1.5e+04 & 1.4e-10 & 4.4e-07 & 1.4e-10	\\
\hline
\end{tabular}
\caption{Numerical results for selected benchmark problems. The table lists the matrix size $n$, preconditioner type, matrix condition number $\kappa(A)$, maximum number of iterations $i$ and selected shift $\sigma^*$ for the shifted p-CG method, and the initial and final residual norms after $i$ iterations for different variants of the CG algorithm.}
\label{tab:benchmarks_specs}
\end{table}


\subsection{Parallel performance results}

This section demonstrates the parallel scalability of the shifted pipelined CG method, and compares to classic CG, pipelined CG and pipelined CG with automated residual replacements (p-CG-rr), which was introduced in \cite{cools2016analysis}. The residual replacement strategy stabilizes the pipelined CG algorithm by incorporating periodic resets of the residual and auxiliary variables to their true values. A small number of additional \textsc{spmv}s is required to compute the corresponding quantities explicitly whenever the residual rounding error becomes too large with respect to the true residual, see \cite{ghysels2014hiding,sleijpen1996reliable,sleijpen2001differences,van2000residual}.

Parallel experiments are performed on a small cluster with $20$ compute nodes, consisting of two $6$-core Intel Xeon 
X5660 Nehalem $2.80$ GHz processors each (12 cores per node). Nodes are connected by $4\,\times\,$QDR 
InfiniBand technology with 32 Gb/s point-to-point bandwidth for message passing and I/O. We use $12$ MPI processes per node to fully exploit parallelism on the machine.
The MPI library used for this experiment is MPICH-3.1.3\footnote{\url{http://www.mpich.org/}}. Note that the environment variables 
\texttt{MPICH\_ASYNC\_PROGRESS=1} and \texttt{MPICH\_MAX\_THREAD\_SAFETY=multiple} are set to ensure optimal parallelism by allowing for non-blocking global communication.
The different variants of the CG algorithm are implemented in PETSc v.3.7.6. The benchmark problem used to asses strong scaling parallel performance is a 2D Poisson model available in PETSc as example $2$ in the KSP folder. The Laplacian operator is discretized using second order finite differences on a $1000\times1000$ grid (1 million DOF). No preconditioner is applied. 

\newpage

\begin{figure}
\begin{center}
\begin{tabular}{cc}
\includegraphics[width=0.45\textwidth]{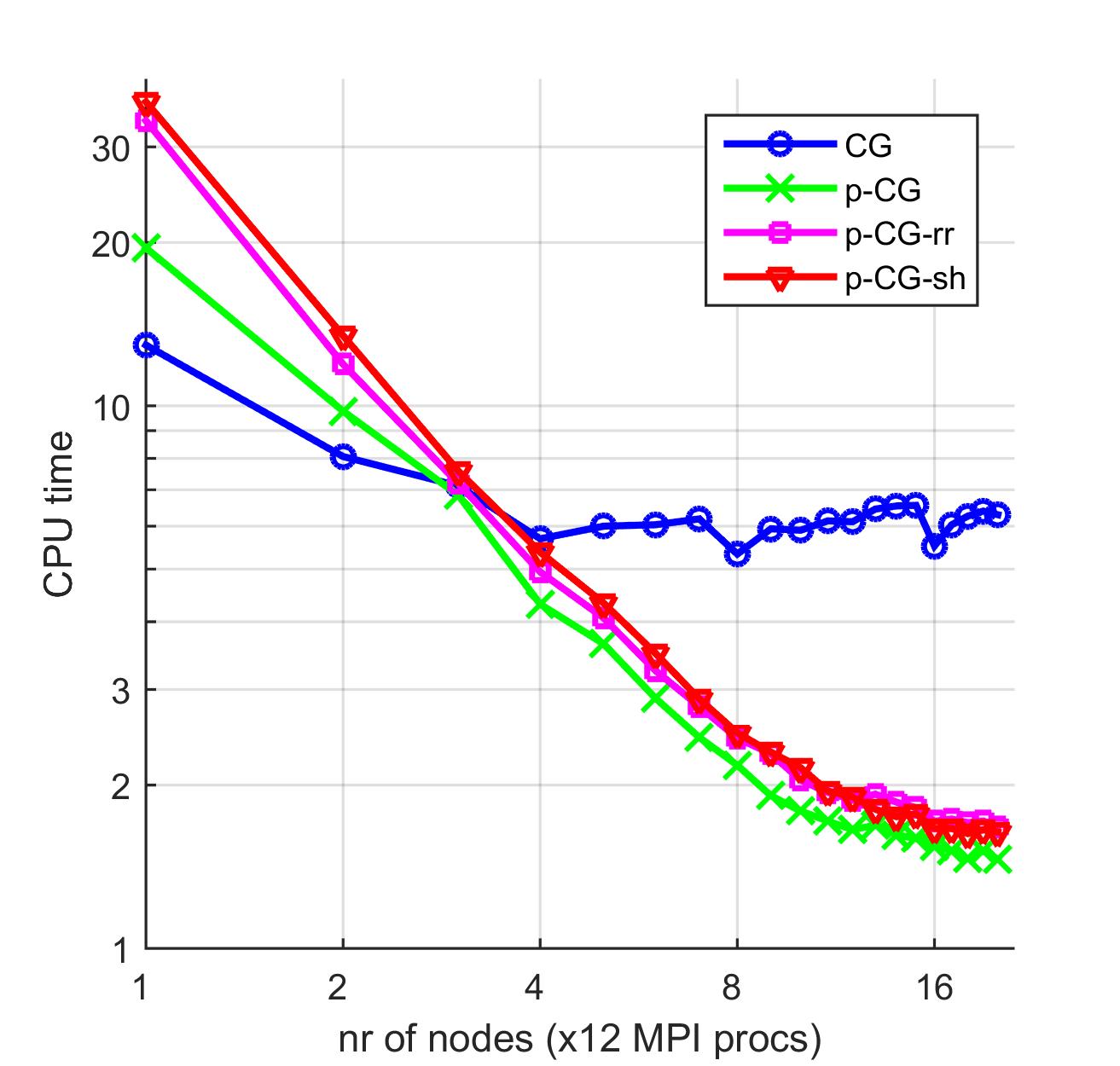} &
\includegraphics[width=0.45\textwidth]{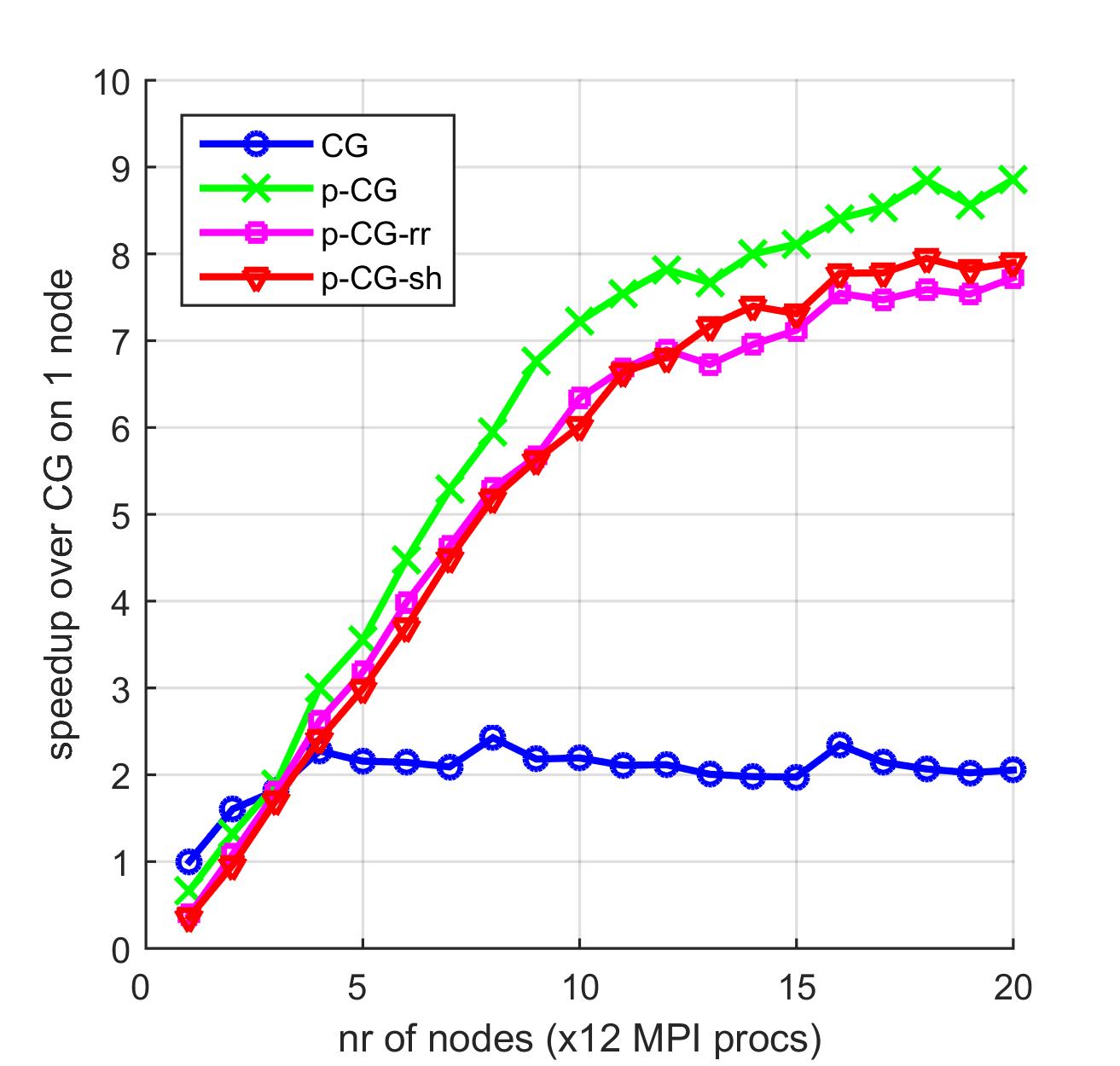} 
\end{tabular}
\end{center}
\caption{Strong scaling experiment on up to $20$ nodes ($240$ cores) for a 2D Poisson problem with $1.000.000$ unknowns.
Left: Absolute time to solution (in seconds) (\texttt{log10} scale) as function of number of nodes (\texttt{log2} scale). 
Right: Speedup over single-node classical CG. 
All methods converged in $1474$ iterations to a relative residual tolerance $1\text{e-}6$; p-CG-rr performed $39$ replacements. 
\label{fig:timings}}
\end{figure}

Figure \ref{fig:timings} shows a strong scaling experiment ranging from 1 to up to 20 nodes. The tolerance imposed on the scaled recursive residual norm $\|\bar{r}_i\|_2 / \|b\|_2$ is $10^{-6}$. The pipelined CG variants clearly out-scale classic CG, achieving a speed-up over single-node CG of approximately $8\times$ on 20 nodes. Classic CG stops scaling at around 4 nodes in this experiment. Performance of the stabilized variants p-CG-rr and p-CG-sh is reduced slightly compared to p-CG due to the additional computations (\textsc{spmv}s for p-CG-rr; \textsc{axpy}s for p-CG-sh) that need to be performed. 

\begin{figure}
\begin{center}
\begin{tabular}{cc}
\includegraphics[width=0.45\textwidth]{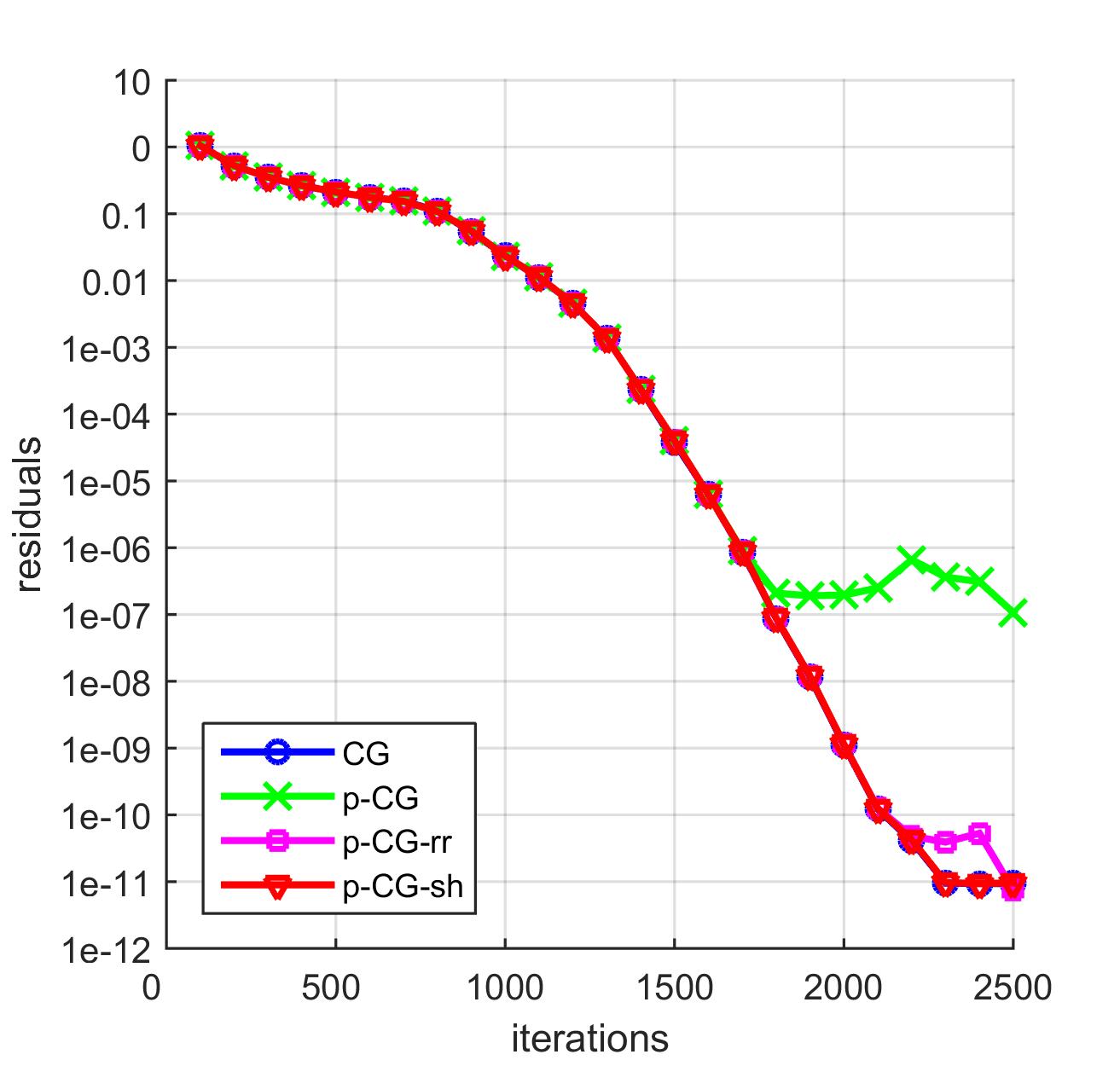} &
\includegraphics[width=0.45\textwidth]{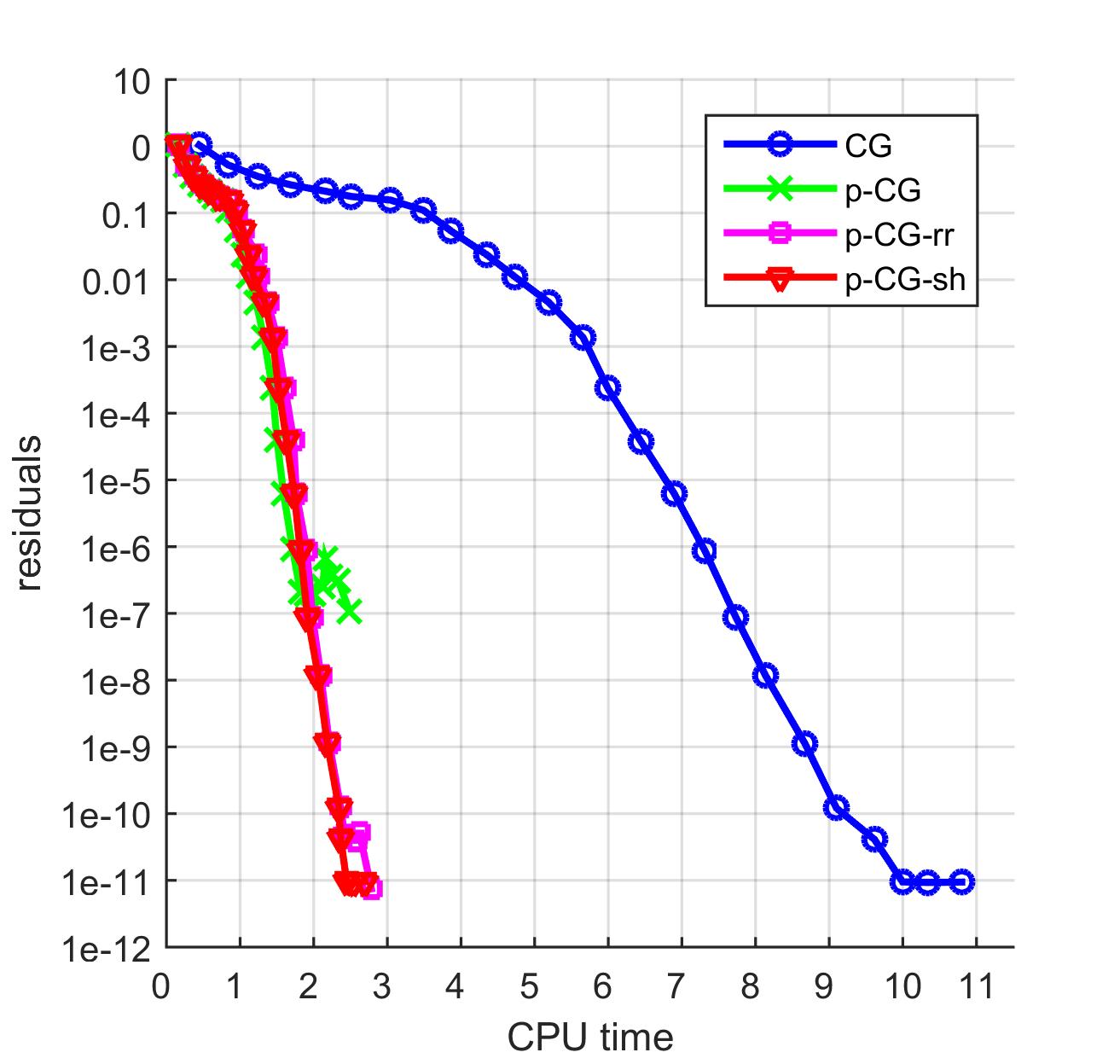} 
\end{tabular}
\end{center}
\caption{Accuracy experiment on $20$ nodes ($240$ cores) for a 2D Poisson problem with $1.000.000$ unknowns.
Left: Explicitly computed residual as function of iterations.
Right: Residual as function of total time spent by the algorithm (in seconds). 
The maximal number of iterations is $2500$ for all methods.
\label{fig:timings2}}
\end{figure}

Figure \ref{fig:timings2} displays accuracy experiments on a 20 node setup. The true residual norm $\|b-A\bar{x}_i\|_2$ is shown in function of the number of iterations (left) and total time to solution (right). Classic CG achieves a high accuracy solution with corresponding residual norm 9.4e-12 in 10.8 seconds. The p-CG method is unable to attain a comparable precision on the final solution, regardless of computational effort. However, it reaches a residual norm around 1.1e-7 in only 2.2 seconds due to the reduction of synchronization bottlenecks and the overlap of global reductions with \textsc{spmv}s. Both p-CG-rr and p-CG-sh are able to attain an accuracy that is comparable to standard CG (with residual norms 7.5e-12 and 9.6e-12 respectively) in around 2.7 seconds, achieving a speedup of approximately $4\times$ over CG while maintaining high accuracy.

\section{Conclusion}

In this research paper we proposed a novel and elegant remedy to the traditionally observed loss of attainable accuracy in the communication hiding pipelined CG algorithm by Ghysels et al.~\cite{ghysels2014hiding}. The proposed methodology is based on a reformulation of the multi-term recurrences for several auxiliary variables in the algorithm. These variables, which exclude the residual, the search direction and the solution itself, are defined using a shifted matrix $(A-\sigma I)$ instead of the original system matrix $A$, and their recurrences are reformulated accordingly. The value of the shift allows to control the build-up of local rounding errors on the solution.

The shifted pipelined CG algorithm is fully equivalent to the classic CG method in exact arithmetic and, with the exception of small number of extra \textsc{axpy}s, does not require any artificial additions to the algorithm. The latter is notable since other common stabilization techniques, such as a residual replacement strategy \cite{cools2016analysis,ghysels2014hiding,sleijpen1996reliable,sleijpen2001differences,van2000residual}, are themselves a possible source of rounding errors, see \cite{liesen2012krylov}. 

The stability analysis presented in this work indicates that the choice of the shifting parameter $\sigma$ is vital for the stability of the pipelined method. Numerical results illustrate that, for a given linear system, a suitable shift can be determined using an \emph{a posteriori} estimate of the local rounding error propagation matrices. The coefficients $\alpha_i$ and $\beta_i$, which are computed as dot-products in each iteration of the (pipelined) CG algorithm, are required to form these propagation matrices. Given a proper choice for the shift parameter, a maximal attainable accuracy on the solution comparable to that of classic CG can be achieved by the shifted pipelined CG method, while parallel scalability is significantly improved compared to CG.

\section{Acknowledgments}

The author acknowledges funding from the Research Foundation Flanders (FWO) under grant application number 12H4617N.

{\footnotesize
\bibliographystyle{unsrt}
\bibliography{refs2}
}

\end{document}